\documentclass{gtart}
\gtart

\usepackage{amsthm}
\usepackage{amsgen}
\usepackage{amsmath,amssymb}
\usepackage{curvesls}
\usepackage{epic}
\usepackage[all]{xy}
\usepackage[dvips]{graphicx}
\usepackage{psfrag}

\newtheorem{thm}{Theorem}[section]

\newtheorem{lem}[thm]{Lemma}
\newtheorem{cor}[thm]{Corollary}
\newtheorem{prop}[thm]{Proposition}
\newtheorem{eg}[thm]{Example}

\theoremstyle{definition}
\newtheorem{defn}[thm]{Definition}

\newcommand{\Complex}{{\mathbb C}}

\newcommand{\I}{{\bf{I}}}
\newcommand{\IG}{{\mathbb G}}

\newcommand{\Real}{{\mathbb R}}
\newcommand{\Natural}{{\mathbb N}}

\newcommand{\Zed}{{\mathbb Z}}

\newcommand{\Diff}{{\mathrm{Diff}}}

\newcommand{\Isom}{{\mathrm{Isom}}}

\newcommand{\splice}{{\bowtie}}

\newcommand{\BAR}[2]{{\lower1.2ex\hbox{#1}\kern-0.5em\text{-----}\kern-0.5em\lower1.2ex\hbox{#2}}}
\newcommand{\Emb}{{\mathrm{Emb}}}

\newcommand{\K}{{\mathcal{K}}}

\newcommand{\RProj}{{\mathbb RP}}

\newcommand{\Brun}{{\mathcal{U}}}
\newcommand{\BrunS}{{\overline{\mathcal{U}}}}

\begin{document}

\title{JSJ-decompositions of knot and link complements in $S^3$}

\author{Ryan Budney}

\address{
Mathematics and Statistics, University of Victoria \\
PO BOX 3045 STN CSC, Victoria, B.C., Canada V8W 3P4 \\
}

\email{rybu@uvic.ca }

\begin{abstract} 
This paper is a survey of some of the most elementary consequences
of the JSJ-decomposition and geometrization for knot and link 
complements in $S^3$.  Formulated in the language of graphs, the result is
the construction of a bijective correspondence between
the isotopy classes of links in $S^3$ and a class of vertex-labelled, 
finite acyclic graphs, called companionship graphs. This construction 
can be thought of as a uniqueness theorem for Schubert's 
`satellite operations.' 
We identify precisely which graphs are companionship graphs of
knots and links respectively.  We also describe how a large family of
operations on knots and links affects companionship graphs. 
This family of operations is called `splicing' and includes, 
among others, the operations of: cabling, connect-sum, 
Whitehead doubling and the deletion of a component.
\end{abstract}


\maketitle

\section{Introduction}\label{INTRODUCTION}

Although the JSJ-decomposition is well-known and frequently used to
study 3-manifolds, it has been less frequently used to study knot complements
in $S^3$, perhaps because in this setting it overlaps with Schubert's `satellite' 
constructions for knots. 
This paper studies the global nature of the JSJ-decomposition for knot
and link complements in $S^3$.
Much of this article is `survey' in nature, in the sense that many of the 
primary results here appear elsewhere in the literature, but not in one place. 
Frequently we offer new proofs of old results, and we attempt to 
refer to the first-known appearance of theorems.

Schubert \cite{Sch2} was the first to study incompressible 
tori in knot complements, which he described in the language of 
`satellite operations.' Among other results, Schubert showed
that satellite operations could be used to recover his connected-sum 
decomposition of knots \cite{Sch1}.
Unlike the case of the connected-sum, 
Schubert did not give a full uniqueness theorem for general 
satellite knots. Waldhausen \cite{Wald_Surf} eventually set up a general
theory of incompressible surfaces in 3-manifolds, which led
to Jaco, Shalen and Johannson's development of the eponymously
named JSJ-decomposition \cite{JacoShalen, Joh2} where an 
appropriate uniqueness theorem was proven.  The JSJ-decomposition 
theorem states that every prime $3$-manifold $M$ contains a collection of embedded, incompressible tori $T \subset M$ so that if one removes an open tubular 
neighbourhood of $T$ from $M$, the resulting manifold $M | T$ is a 
disjoint union of Seifert-fibred and atoroidal manifolds.  
Moreover, if one takes a minimal collection of such tori, they are 
unique up to isotopy.  It is the purpose of this paper to work out 
the explicit consequences of this theorem and the later developments in 
Geometrization, for knot and link complements in $S^3$. 

Given a link $L \subset S^3$, with complement $C_L= S^3 \setminus U$ 
(where $U$ is an open tubular neighbourhood of $L$) a collection of natural 
questions one might ask is:

\begin{enumerate}
\item Which Seifert-fibred manifolds can be realised as components of
$C_L | T$ for $T$ the JSJ-decomposition of a knot or link complement $C_L$?
\item Which non Seifert-fibred manifolds arise in the same way?
\item How are the above manifolds embedded in $S^3$?
\item How do they all `fit together' globally, and what combinations are
possible?
\end{enumerate}

We partially answer item (3) first. We prove in Proposition \ref{embtorbou} 
that if $M$ is a compact submanifold of $S^3$ with $\partial M$ a disjoint union of $n$ embedded tori, 
if we let $p$ and $q$ be the number of solid tori components and non-trivial
knot complement components of $S^3 \setminus int(M)$ respectively, 
where $p+q=n$ then there exists an embedding $f:M \to S^3$ so that 
$f(M)$ is the complement of an open tubular neighbourhood of an $n$-component link  
$L \subset S^3$ which contains a $q$-component unlink as a sublink. This
brings Brunnian properties into the picture. 
We go on to prove in Proposition \ref{desplice_conv} that 
there is a canonical choice for $f$, thus the study of
submanifolds of $S^3$ with torus boundary reduces in a natural
way to link theory in $S^3$.

Section \ref{SEIFERT} answers question (1) by computing which Seifert-fibred
manifolds embed in $S^3$ in Proposition \ref{sf2}. 
This allows us, in Proposition \ref{seifert-glt} to determine the links  in 
$S^3$ which have Seifert-fibred complements. We develop a notation
sufficient to describe all links with Seifert-fibred complements,
and classify them up to unoriented isotopy in Proposition \ref{Spquniq}.
Proposition \ref{FSEP} gives conditions on when two
Seifert-fibred manifolds can be adjacent in the JSJ-decomposition of
a link complement, and we end Section \ref{SEIFERT} with a discussion
of the geometric structures on Seifert-fibred link complements. 

Section \ref{MAIN_RESULTS} is the heart of the paper where we investigate item (4).
The components of $C_L | T$ naturally form the vertex-set
of a partially-directed acyclic graph $G_L$, called the JSJ-graph of 
$L$ (Definition \ref{def1}). 
In Definition \ref{COMP_GRAPH_DEF} we construct the companionship graph 
$\IG_L$ of a link $L$ by labelling the
vertices of $G_L$ with natural `companion links' to $L$.
We show $\IG_L$ is a complete isotopy invariant of $L$ in
Proposition \ref{IG_L_is_complete}. 
From here we investigate the properties of the graphs $\IG_L$. 
The most basic property of $\IG_L$ is that it is a `splicing
diagram' (Definition \ref{splice_diagram}).  Roughly, this means
that $\IG_L$ is a finite, acyclic graph, some of whose edges are
oriented, where each vertex labelled by a link, and each edge of
the graph corresponds to a matched pair of components of the links
decorating the endpoints of the edge.
We complete the study of companionship graphs of knots first,
giving a characterisation of the companionship graphs of knots
in Theorem \ref{mainthm}, showing, among other things, that
the companionship graphs are naturally rooted trees.  
To characterise companionship trees of knots, we use a notion of splicing
(Definition \ref{splicedef}) that allows us to inductively construct
knots with ever more complicated companionship trees.
Returning to links, we show that $\IG_L$ satisfies a local
Brunnian property (Proposition \ref{brun_excl}).  Splicing diagrams
that satisfy this property we call valid. We show in Proposition
\ref{REAL_VALID} that valid splicing diagrams correspond bijectively
to isotopy classes of collections of disjoint embedded tori in link complements. This allows us to identify in Proposition \ref{graph_real}
which valid splice diagrams are companionship graphs.
Similarly to the case of links, we use splicing to inductively
construct companionship graphs. In Proposition \ref{graph_of_a_splice}
we show how splicing affects the companionship graphs of an
arbitrary pair of links.

Thurston's Hyperbolisation Theorem answers question (2), telling 
us that the interiors of non-Seifert-fibred components of 
$C_L | T$ are complete hyperbolic manifolds of finite-volume.
At the end of Section \ref{MAIN_RESULTS} we describe algorithms for finding
the JSJ-decomposition and the geometric structures on the components.

Our definitions and conventions regarding knots and links are:
\begin{itemize}
\item A knot is a compact, connected, boundaryless, oriented, $1$-dimensional sub-manifold of $S^3$. Thus a knot is diffeomorphic to $S^1$.
\item Given a knot $K$ let $+K=K$ and let $-K$ be the oppositely-oriented
knot, meaning as unoriented manifolds $-K$ and $K$ are the same, thus $-(-K)=K$
but $-K \neq K$.
\item Given a finite set $A$, a link $L$, indexed by $A$, is
a disjoint collection of knots $\{L_a : a \in A\}$.  
Given $A' \subset A$, $L_{A'}$
denotes the sublink of $L$ indexed by $A'$.
\item Given two links $L$ and $L'$ with index-set $A$,
an isotopy from $L$ to $L'$ is an orientation-preserving
diffeomorphism $f : S^3 \to S^3$ such that $f(L_a)=L'_a$ for all
$a \in A$. This notion agrees with the traditional notion of isotopy
because all orientation-preserving diffeomorphisms of $S^3$ are isotopic
to the identity by Cerf \cite{Cerf}.
\item Given two links $L$ and $L'$ with index sets $A$ and $A'$,
we say $L$ and $L'$ are unoriented-isotopic if they are isotopic as unoriented
submanifolds of $S^3$. Stated another way, an orientation-preserving 
diffeomorphism $f : S^3 \to S^3$ is called
an unoriented isotopy from $L$ to $L'$ if there is
a bijection $\sigma : A \to A'$ and a function $\epsilon_f : A \to \{\pm\}$
such $f(L_a) = \epsilon_f(a) L'_{\sigma(a)}$ for all $a \in A$.
Provided $\epsilon_f(a)=+$ for all $a \in A$, we say $L$ and $L'$ are
isotopic modulo re-indexing.
\item Let $D^n = \{ x \in \Real^n : |x| \leq 1 \}$ is the compact unit $n$-disc.  
\item Given a solid torus $M \simeq S^1 \times D^2$ in $S^3$, there are two
canonical isotopy classes of unoriented curves in $\partial M$, the meridian and 
longitude respectively.  The meridian is the essential curve in $\partial M$ that 
bounds a disc in $M$.  The longitude is the essential curve in $\partial M$ that 
bounds a $2$-sided surface in $S^3 \setminus int(M)$ (a Seifert surface). 
If $M$ is a closed tubular neighbourhood of a knot $K$, the longitude of
$M$ is parallel to $K$ thus we give it the induced orientation. We give
the meridian $m$ the orientation so that $lk(K,m)=+1$.
\item For such standard definitions as connected-sum, splittability, etc, 
we will follow the notation of Kawauchi \cite{Kaw}.
\item Following the conventions of Kanenobu \cite{Kan} and Debrunner \cite{DeBrun}, given $L$ 
indexed by $A$ we define $\Brun_L \subset 2^A$ by the rule 
$\Brun_L = \{ S \subset A : L_S \text{ is not split}\}$. 
$\Brun_L$ is the Brunnian property of $L$.  
\item $\BrunS_L = \{ S \subset A : L_S \text{ is a } |S|\text{-component unlink}\}$
is called the strong Brunnian property of $L$.
\end{itemize}

Our definitions and conventions regarding $3$-manifolds are:
\begin{itemize}
\item $3$-manifolds are taken to be oriented and are allowed boundary. 
\item For standard definitions of connected-sum, prime, irreducible, Seifert-fibred,
incompressible surface, etc, we will use the conventions of Hatcher \cite{Hatcher3}.   
\item Given a $3$-manifold $M$ and a properly-embedded $2$-sided surface 
$S \subset M$ define $M | S = \{V \subset M : V=\overline{W} 
\text{ where } W \text{ is a path-component of } M \setminus S\}$. 
We call the elements of $M | S$ `the components of $M | S$' 
or `the components of $M$ split along $S$'.
If $S'$ is a component of $S$, then $S'$ is a submanifold of at most
two components of $M|S$.
\item Given $T$ a collection of disjoint embedded incompressible tori 
in an irreducible $3$-manifold $M$ we say $T$ is the 
$JSJ$-decomposition if $M | T$ is a disjoint
union of Seifert-fibred and atoroidal manifolds, 
and if no smaller such collection of tori splits $M$ into Seifert-fibred
and atoroidal manifolds. 
\end{itemize}

There are several treatments of the JSJ-decompositions of $3$-manifolds. 
There's the original work of Jaco and Shalen \cite{JacoShalen} and the
simultaneous work of Johannson \cite{Joh2}.  Some more recent
expositions are available as well. There is Hatcher's notes 
\cite{Hatcher3}, which this article follows, and also 
the notes of Neumann and Swarup \cite{Neumann}.  
As Neumann and Swarup \cite{Neumann} point out (Proposition 4.1) 
all versions of the decomposition are closely related. 
In the case of $3$-manifolds $M$ with $\chi (M) = 0$ 
such as link complements in $S^3$, all versions of the 
JSJ-decomposition are the same.  When $\chi (M) \neq 0$
one can get various different incompressible annuli in the
decompositions, depending on whose conventions are followed.
This difference is important as the original 
JSJ-decomposition is more closely related to Thurston's 
geometric decomposition of 3-manifolds.

There are also treatments of the JSJ-decomposition of knot and link 
complements.  The book of Eisenbud and Neumann \cite{EN} gives a
detailed analysis of the structure of the JSJ-decomposition of 
links in homology spheres whose complements are graph manifolds.  
Our paper differs from their book in that we study the class of
links in $S^3$ with no restriction on the complements.
The aspect of this paper which is new is that the complicating
factor is the Brunnian properties of the resulting companion links. 
A once frequently-quoted yet unpublished manuscript of
Bonahon and Siebenmann's \cite{BS} also investigates JSJ-decompositions
of link complements in $S^3$.  Various results from the manuscript of Bonahon
and Siebenmann appears in the survey of Kawauchi's \cite{Kaw}, 
but the results that do appear are on the $\Zed_2$-equivariant 
JSJ-decomposition of the $2$-sheeted branched cover of $S^3$ over a link, which 
they show correspond to Conway Spheres. The current exposition is most
closely related to the part of the Bonahon-Siebenmann manuscript \cite{BS} that
does not appear in Kawauchi's survey.
This paper in part duplicates the results
of Schubert, Eisenbud, Neumann, Bonahon and Siebenmann, and we attempt
to give full credit to their discoveries.

This article started out as a technical lemma needed to determine the
class of hyperbolic $3$-manifolds that appear as components of a knot
complement split along its JSJ-tori.  The observation that the result,
although known to some, is not `well known,' motivated the author to
put together the present exposition.
I would like to thank Allen Hatcher for several early suggestions on how
to approach the topic. I'd also like to thank Gregor Masbaum for his
suggestion of reformulating what is now Proposition \ref{unknottedcomp}, which
led to the connection with Schubert's paper \cite{Sch2}. I'd like to thank Daniel
Moskovich and David Cimasoni for their comments on the paper, and 
Francis Bonahon for encouragement.

\section{Disjoint knot complements, and companions}\label{COMPLEMENTS}

This section starts with a technical proposition about disjoint knot complements in $S^3$
which ultimately motivates the notions of splicing and companions of a knot or link.

\begin{prop}\label{unknottedcomp} Let $C_1, C_2, \cdots, C_n$ be $n$ disjoint submanifolds of $S^3$
such that for all $i \in \{1,2,\cdots,n\}$, $K_i = \overline{S^3 \setminus C_i}$ is a non-trivially embedded solid-torus in $S^3$.  Then there exists $n$ disjointly embedded
$3$-balls $B_1,B_2,\cdots,B_n \subset S^3$ such that $C_i \subset B_i$ for all
$i \in \{1,2,\cdots,n\}$. Moreover, each $B_i$ can be chosen to be $C_i$ union a
$2$-handle which is a tubular neighbourhood of a meridional disc for $K_i$.
\begin{proof}
For all $i \in \{1,2,\cdots,n\}$ let
$D^2_i$ be a meridional disc for $K_i$. 

Consider the case that we have $j$ disjoint $3$-balls $B_1, B_2, \cdots, B_j$
such that $C_i \subset B_i$ for all $i \in \{1,2,\cdots,j\}$ with $B_i$ 
disjoint from $C_l$ for all $i \neq l$, $i \in \{1,2,\cdots,j\}$ and $l \in \{1,2,\cdots,n\}$. 
We proceed by induction, the base-case being the trivial $j=0$ case.

Consider the intersection of $D^2_{j+1}$
with $\partial B_1 \cup \partial B_2 \cup \cdots \cup \partial B_j \cup 
\partial C_{j+2} \cup \partial C_{j+3} \cup \cdots \cup \partial C_n$.

\begin{itemize}
\item If the intersection is empty, let $B_{j+1}$ be a
regular neighbourhood of $C_{j+1} \cup D^2_{j+1}$. 
\item If on the other hand the intersection is non-empty, let $S$ be an innermost
curve of the intersection bounding an innermost disc $D$ in $D^2_{j+1}$. 
Thus $S$ is a sub-manifold
of one of $\partial B_1, \cdots, \partial B_j$ or 
$\partial C_{j+2}, \cdots, \partial C_n$. 
\begin{itemize}
\item If $S$ bounds a disc $D'$ in some $\partial B_i$ for 
$i \leq j$ or $\partial C_i$ for
$j+2 \leq i \leq n$, then $D \cup D'$ bounds a ball in $B_i$ or $C_i$ respectively,
which gives a natural isotopy of $D^2_{j+1}$ which lowers
the number of components of intersection with the family  
$\partial B_1 \cup \partial B_2 \cup \cdots \cup 
\partial B_j \cup \partial C_{j+2} \cup \partial 
C_{j+3} \cup \cdots \cup \partial C_n$.
\item If $S$ does not bound a disc in the above family, it must be
a meridional curve in some $\partial C_i$ for $j+2 \leq i \leq n$. In this case, 
we let $B_i$ be a regular neighbourhood of $C_i \cup D$. 
\end{itemize} 
\end{itemize}
Thus by re-labelling the tori and balls appropriately, we have completed the
inductive step.
\end{proof}
\end{prop}

Proposition \ref{unknottedcomp} first appears in the literature as a theorem of
Schubert's \cite{Sch2} (\S 15.1, pg. 199). It also appears in the work
of Sakuma's \cite{Sakuma} on the symmetry properties of knots. 
A related result was re-discovered by Bonahon and Siebenmann in their 
unpublished manuscript \cite{BS} as part of their
algorithm to determine if a `splicing tree' results in the construction of a link
embedded in $S^3$. 

\begin{prop}\label{embtorbou}
Let $M$ be a compact submanifold of $S^3$ with 
$\partial M$ a disjoint union of $n$ 
tori. By Alexander's Theorem, 
$\overline{S^3 \setminus M}$ consists of a 
disjoint union of $p$ solid tori and $q$ non-trivial knot 
complements, where $p+q=n$. 
There exists an embedding $f:M \to S^3$ such
that $f(M)$ is the complement
of an open tubular neighbourhood of an $n$-component 
link $L \subset S^3$ which contains a 
$p$-component unlink as a sublink.
\begin{proof}
Let $q \in \{0, 1, \cdots, n\}$ be the number of components 
of $\overline{S^3 \setminus M}$ which are non-trivial knot
complements, and let
$\overline{S^3 \setminus M} = C_1 \sqcup \cdots \sqcup C_n$ where 
$C_i$ is a solid torus for $q+1 \leq i \leq n$ and
a non-trivial knot-complement for $1 \leq i \leq q$.
By Proposition \ref{unknottedcomp} there exists disjoint 
$3$-balls $B_1, \cdots, B_q \subset S^3$ such that $B_i$ is obtained from
$C_i$ by an embedded $2$-handle attachment, $B_i = C_i \cup H^2_i$. Dually,
$C_i$ is obtained from $B_i$ by drilling out a neighbourhood of a 
knotted properly-embedded interval. 

Let $Q = \{(x,y,z) \in \Real^3 : |(y,z)| \leq \frac{1}{2} \} \cap D^3$.
For $i \in \{1,2,\cdots,q\}$ let 
$u_i : (D^3,A_i) \to (B_i,H_i^2)$ be an orientation-preserving
diffeomorphism of pairs,
where $A_i \subset D^3$ is a $3$-ball such that 
$A_i \cap \partial D^3 = Q \cap \partial D^3$.
For $i \in \{1,2,\cdots,q\}$ let 
$w_i : Q \to A_i$ be a diffeomorphism which
is the identity on $Q \cap \partial D^3 = A_i \cap \partial D^3$. 

We define an embedding $f : M \to S^3$ as follows:

$$f(x)= \left\{
\begin{array}{lr}
x & \text{ if } x \in \overline{M \setminus \cup_{i=1}^q H^2_i} \\
u_i \circ w_i^{-1}\circ u_i^{-1}(x) & \text{ if } x \in H^2_i \\
\end{array}
\right.$$

By design, $f(\cup_{i=1}^q \partial C_i)$ bounds a tubular neighbourhood of
a $q$-component unlink in the complement of $f(M)$.
To argue that $f(M)$ is a link complement, notice that in our
definition $f$ extends naturally to an embedding 
$\overline{S^3 \setminus \cup_{i=1}^q C_i} \to S^3$.
\end{proof}
\end{prop}

\begin{defn}\label{defperiph}
Let $M \subset S^3$ be a $3$-manifold, and let $T \subset \partial M$
be a torus. Provided $C$ is the component of $\overline{S^3 \setminus M}$ containing $T$,
an essential curve $c \subset T$ is called an external (resp. internal) 
peripheral curve for $M$ at $T$ if $c = \partial S$ for some properly-embedded 
surface $S \subset C$ (resp. $S \subset \overline{S^3 \setminus C}$).
\end{defn}

\begin{prop}\label{desplice_conv}
Let $M \subset S^3$ be a $3$-manifold whose boundary is a disjoint
union of tori.  Up to isotopy, there exists a unique orientation-preserving
embedding $f : M \to S^3$ such that (1) and (2) are true:
\begin{enumerate}
\item $f(M)$ is the complement of a tubular neighbourhood of a link in $S^3$.
\item $f$ maps external peripheral curves of $\partial M$ to external peripheral
curves of $\partial(f(M))$.
\end{enumerate}
$f$ will be called the untwisted re-embedding of $M$.
\begin{proof}
We prove existence in the framework of the proof of Proposition \ref{embtorbou}.
Let $Q' =\{(x,y,z) \in \partial D^3 : -\frac{\sqrt{3}}{4} \leq x \leq \frac{\sqrt{3}}{4} \}$. 
For each $i \in \{1,2,\cdots,q\}$ extend $w_i$ to be
a homeomorphism $w_i : Q\cup Q' \to A_i \cup Q'$ so that $w_i(x) =x$ for
all $x \in Q'$.
Fix the curves 
$c_1 = \{ (\cos \theta,\sin \theta,0) : \pi \leq \theta \leq 2\pi \} \cup
\{(x,0,0) : -1 \leq x \leq 1\}$ and 
$c_2 = \{ (\cos \theta,\sin \theta,0) : \frac{\pi}{6} \leq \theta \leq \frac{5\pi}{6} \}
\cup \{(x,\frac{1}{2},0) : -\frac{\sqrt{3}}{4} \leq x \leq \frac{\sqrt{3}}{4}\}$.
For each $i \in \{1,2,\cdots,q\}$ there is a unique choice of $w_i$ up to 
isotopy so that $lk(w_i(c_1),w_i(c_2))=0$, since any two choices of $w_i$ differ
by some Dehn twist about a disc in $Q$ separating the components of
$Q \cap \partial D^3$. With this choice, then external peripheral curves are sent to external peripheral curves.

To prove uniqueness, let $f_1 : M \to S^3$ and $f_2 : M \to S^3$ be any
two embeddings of $M$ in $S^3$ as link complements sending external 
peripheral curves to external peripheral curves. Let
$\overline{M}$ be the Dehn filling of $M$ where the attaching maps are
given by the external peripheral curves of $M$. Then $f_1$ and $f_2$ extend
to orientation-preserving diffeomorphisms 
$\overline{f_1}, \overline{f_2} : \overline{M} \to S^3$. Cerf's
theorem \cite{Cerf} states that any two orientation-preserving 
diffeomorphisms of $S^3$ are isotopic, thus $f_1$ and $f_2$ are isotopic.
\end{proof}
\end{prop}

\begin{defn}\label{defcompanionlink}
Given $L$ a link in $S^3$ let $C_L = S^3 \setminus U_L$ 
where $U_L$ is any open tubular neighbourhood of $L$.
If $M \subset C_L$ is a manifold with incompressible torus boundary
let $f : M \to S^3$ be its untwisted
re-embedding. Then $f(M)$ is the complement of some link $L'$ in $S^3$,
unique up to unoriented isotopy. Any such link will be
called a companion link to $L$. If $M$ is a component of one of
the prime summands of $C_L$
split along the tori of its JSJ-decomposition, we call
$L'$ a JSJ-companion link to $L$.
A link $L$ will be called compound if $C_L$ is reducible
or if it is irreducible with a non-empty JSJ-decomposition.
\end{defn}

Thus a link is non-compound if and only if it is atoroidal and not
split.

\section{Seifert-fibred submanifolds of $S^3$}\label{SEIFERT}

In this section we determine which Seifert-fibred manifolds embed
in $S^3$, and the various ways in which they embed.  This allows us to 
classify the links in $S^3$ whose complements are Seifert-fibred,
and give basic restrictions on which Seifert-fibred manifolds can be
adjacent in the JSJ-decomposition of a $3$-manifold in $S^3$.

\begin{lem}\label{sfl} If $M$ be a sub-manifold of $S^3$ with
non-empty boundary a union of tori, 
then either $M$ is a solid torus or a 
component of $\overline{S^3 \setminus M}$ is a solid torus.
\begin{proof}
Let $C=\overline{S^3 \setminus M}$. Since $\partial M$ consists of a disjoint union
of tori, every component of $\partial M$ contains an essential curve $\alpha$ which bounds a disc $D$ in $S^3$.  Isotope $D$ so that it intersects $\partial M$ transversely
in essential curves. Then $\partial M \cap D \subset D$ consists of a finite collection
of circles, and these circles bound a nested collection of discs in $D$. Take an innermost
disc $D'$. If $D' \subset M$ then $M$ is a solid torus. If $D' \subset C$ then
the component of $C$ containing $D'$ is a solid torus.
\end{proof}
\end{lem}

We will use the notation in Hatcher's notes \cite{Hatcher3} for
describing orientable Seifert-fibred manifolds. In short, 
let $\mathcal{S}_{g,b}$ denote a compact 
surface of genus $g$ with $b$ boundary components.
If $g<0$, $\mathcal{S}_{g,b}$ is 
the connect-sum of $-g$ copies of $\RProj^2$ and $\mathcal{S}_{0,b}$.
$M(g,b;\frac{\alpha_1}{\beta_1},\frac{\alpha_2}{\beta_2},
\cdots,\frac{\alpha_k}{\beta_k})$
denotes the Seifert-fibred manifold fibred over $\mathcal{S}_{g,b}$
with at most $k$ singular fibres, and fibre-data $\alpha_i/\beta_i$.
One constructs $M(g,b;\frac{\alpha_1}{\beta_1},\frac{\alpha_2}{\beta_2},
\cdots,\frac{\alpha_k}{\beta_k})$ from the orientable
$S^1$-bundle over $\mathcal{S}_{g,b+k}$ by Dehn filling along $k$ of 
the boundary components using the attaching slopes 
$\frac{\alpha_i}{\beta_i}$ for $i \in \{1,2,\cdots,k\}$.

%\begin{defn}\label{standarddef} 
We give a non-standard but flexible notation for defining unoriented 
isotopy classes of links in $S^3$
which are the union of fibres from a Seifert fibring
of $S^3$. Provided $(p,q)\in \Zed^2$ satisfies $p \neq 0$ and
$q \neq 0$, $S(p,q|X)$ denotes the subspace of $S^3$ made up of the union of
three disjoint sets $S_1$, $S_2$, $S_3$ where:
\begin{itemize}
\item $S_1 = \{(z_1,z_2) \in S^3 \subset \Complex^2 : z_1^{p}=z_2^{q}\}$.
\item $S_2 = \{(z_1,0) \in S^3 \}$ provided $*_1 \in X$, otherwise
$S_2 = \emptyset$. 
\item $S_3 = \{(0,z_2) \in S^3 \}$ provided $*_2 \in X$,
otherwise $S_3=\emptyset$. 
\end{itemize}
%\end{defn}

We mention some shorthand notation for some common links with 
Seifert-fibred complements. 
The Hopf link $H^1$ is the $2$-component link in $S^3$ given by
the union of $S_2$ and $S_3$ above, alternatively
$H^1=S(2,2|\emptyset)$. Up to isotopy (modulo re-labeling) 
there are two Hopf links,
distinguished by the linking number of the components.
{\psfrag{hopf}[tl][tl][1][0]{$H^1$}
$$\includegraphics[width=4cm]{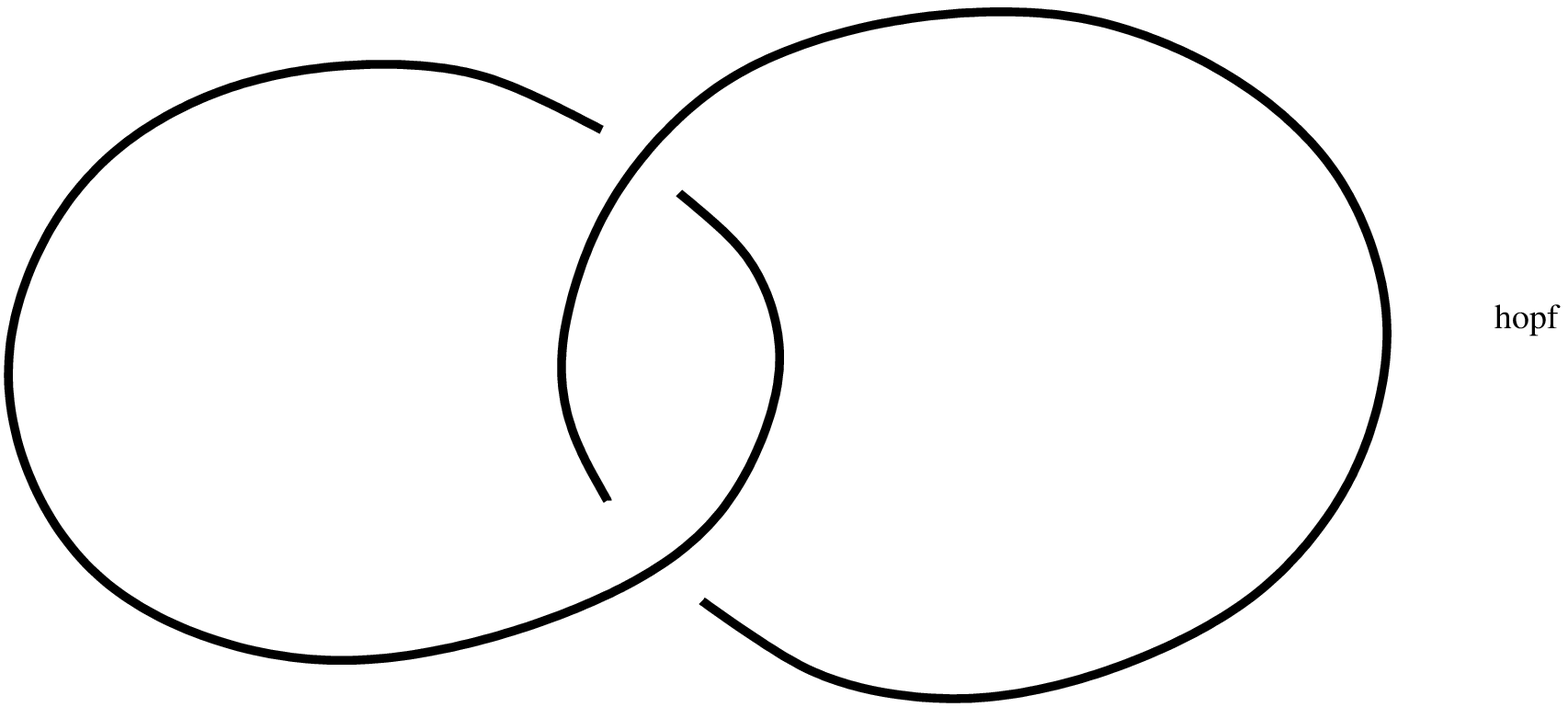}$$}
If one takes a connected-sum of $p$ copies of the Hopf link, one obtains
the $(p+1)$-component `key-chain link' $H^p$. Module re-labelling, 
there are $p+1$ $(p+1)$-component key-chain links.  When
orientation matters, we will use the `right handed' key-chain link
where all non-zero linking numbers are positive.
{\psfrag{...}[tl][tl][0.6][0]{$\cdots$}
\psfrag{hopf}[tl][tl][1][0]{$H^p$}
$$\includegraphics[width=4cm]{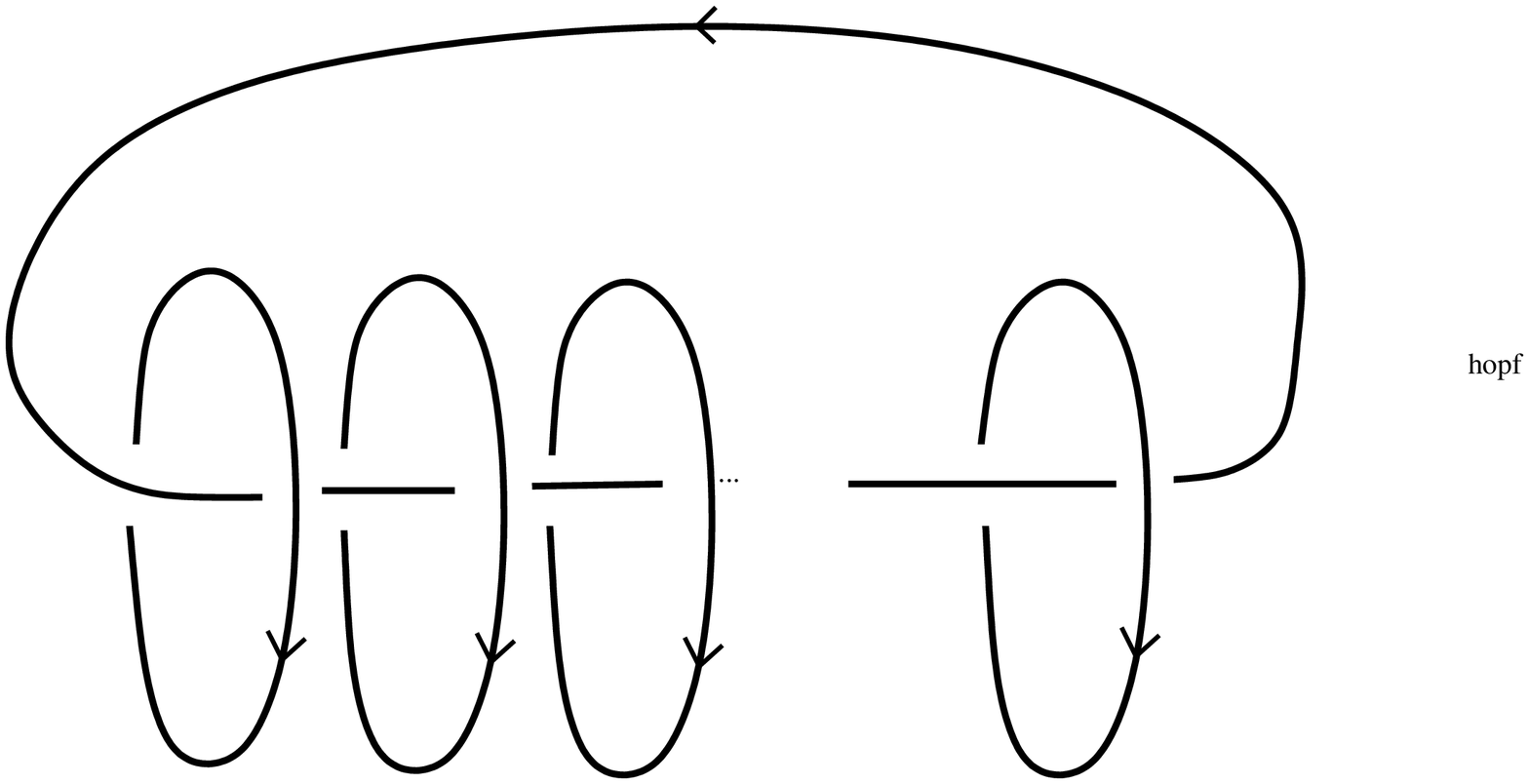}$$}
For any $(p,q) \in \Zed^2$ with
$p \nmid q$, $q \nmid p$ and $GCD(p,q)=1$ the $(p,q)$-torus knot
is $T^{(p,q)} = S(p,q|\emptyset)$. There is
only one $(p,q)$-torus knot, since all torus knots are are invertible.
The conditions $q \nmid p$ and $p \nmid q$ ensure the unknot is
not a torus knot.

For any $(p,q) \in \Zed \times \Zed$ with $p\nmid q$ and $GCD(p,q)=1$,
the $(p,q)$-Seifert link $S^{(p,q)}$ is defined to be
$S^{(p,q)} = S(p,q|\{*_1\})$. We fix the orientation on $*_1$
counter-clockwise and orient the remaining component by the parametrisation
$(\frac{z^q}{\sqrt{2}},\frac{z^p}{\sqrt{2}})$ where $z \in S^1$.
Our condition $p \nmid q$ is there to ensure that the Hopf link is not
considered to be a Seifert link.
{
\psfrag{z1}[tl][tl][1][0]{$z_1$}
\psfrag{z2}[tl][tl][1][0]{$z_2$}
\psfrag{spq}[tl][tl][1][0]{$S^{(p,q)}$}
\psfrag{s1}[tl][tl][1][0]{$*_1$}
$$\includegraphics[width=6cm]{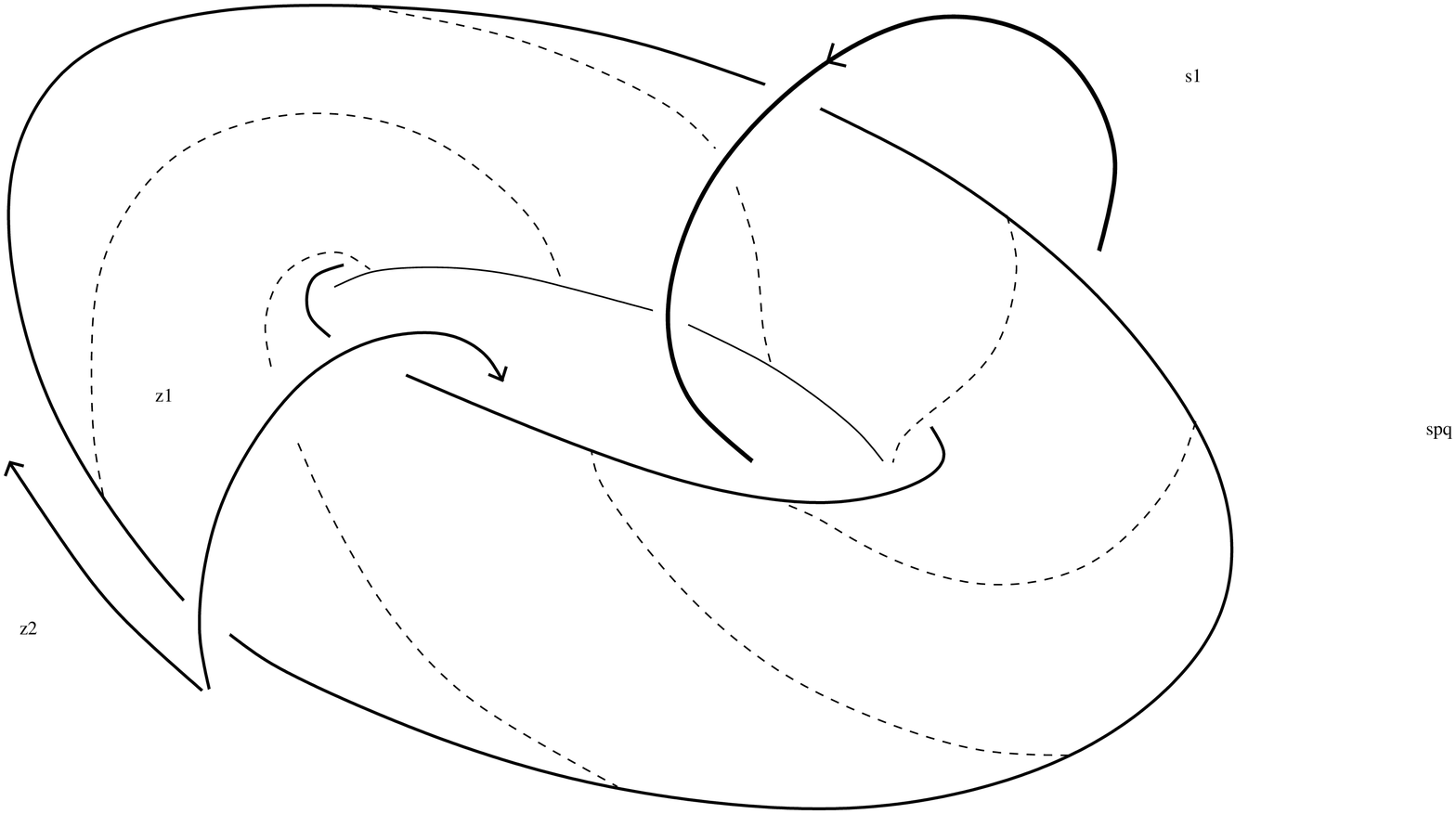}$$
}
\begin{prop}\label{sf2} Let $V \neq S^3$ be a Seifert-fibred sub-manifold of $S^3$, 
then $V$ is diffeomorphic to one of the following:
\begin{itemize}
\item $M(0,n;\frac{\alpha_1}{\beta_1},\frac{\alpha_2}{\beta_2})$ for $n \geq 1$ and
$\alpha_1\beta_2-\alpha_2\beta_1=\pm 1$. These appear only as the complements of $n$ 
regular fibres in a Seifert fibring of $S^3$.  
\item $M(0,n;\frac{\alpha_1}{\beta_1})$ for $n \geq 1$. These appear 
only as complements of $n-1$ regular fibres in a Seifert-fibring of an 
embedded solid torus in $S^3$.
\item $M(0,n;)$ for $n \geq 2$. These appear in two different ways:
\begin{itemize}
\item As complements of the singular fibre and $n-1$ regular fibres 
in a Seifert-fibring an embedded solid torus in $S^3$.
\item A manifold whose untwisted re-embedding in $S^3$ is the complement of 
a key-chain link. 
\end{itemize}
\end{itemize}
\begin{proof}
Consider $V \simeq M(g,b;\frac{\alpha_1}{\beta_1},\frac{\alpha_2}{\beta_2},
\cdots,\frac{\alpha_k}{\beta_k})$

By design $b=0$ if and only if $V=S^3\simeq M(0,0;\frac{\alpha_1}{\beta_1},\frac{\alpha_2}{\beta_2})$ where $\alpha_1\beta_2-\alpha_2\beta_1 = \pm 1$.

Consider the case $b \geq 1$.
Seifert-fibred manifolds that fibre over a non-orientable
surface do not embed in $S^3$ since an orientation-reversing closed
curve in the base lifts to a Klein bottle, which does not embed in $S^3$
by the Generalised Jordan Curve Theorem \cite{Pollack}. Thus $g \geq 0$. 

A Seifert-fibred manifold that fibres over a surface of genus $g>0$ does
not embed in $S^3$ since the base manifold contains two curves that
intersect transversely at a point. If we lift one of these curves
to a torus in $S^3$, it must be non-separating. This again contradicts
the Generalised Jordan Curve Theorem, thus $g=0$. 

By Lemma \ref{sfl}, either $V$ is a solid torus $V \simeq M(0,1;\frac{\alpha_1}{\beta_1})$
or some component $Y$ of $\overline{S^3 \setminus V}$ is a solid torus. Consider 
the latter case.  There are two possibilities. 

\begin{enumerate}
\item The meridians of $Y$ are fibres of $V$. 
If there is a singular fibre in $V$, let $\beta$ be an embedded arc in the
base surface associated to the Seifert-fibring of $V$ which starts at the singular
point in the base and ends at the boundary component corresponding to $\partial Y$.
$\beta$ lifts to a $2$-dimensional CW-complex in $M$, and the endpoint of
$\beta$ lifts to a meridian of $Y$, thus it bounds a disc.  If we append this
disc to the lift of $\beta$, we get a CW-complex $X$ which consists of a $2$-disc
attached to a circle. The attaching map for the $2$-cell is multiplication by
$\beta$ where $\frac{\alpha}{\beta}$ is the slope associated to the singular fibre. 
The boundary of a regular neighbourhood of $X$ is a $2$-sphere, so we have 
decomposed $S^3$ into a connected sum $S^3 = L_{\frac{\gamma}{\beta}} \# Z$ 
where $L_{\frac{\gamma}{\beta}}$ is a lens space with $H_1 L_{\frac{\gamma}{\beta}} = \Zed_\beta$. Since $S^3$ is irreducible, $\beta=1$.
Thus $V\simeq M(0,n;)$ for some $n \geq 1$.  Consider the untwisted
re-embedding $f : V \to S^3$, and let $V' = f(V)$.
$V'$ is the complement of some link $L=(L_0,L_1, \cdots,L_{n-1})$
such that $(L_1,\cdots,L_{n-1})$ is an unlink 
(provided we let $L_0$ correspond to $Y$).
$C_{L_0}$ is obtained from $V'$ by Dehn filling on $L_1, \cdots, L_n$
with integral slopes, thus $C_{L_0}$ is a solid torus.
\item The meridians of $Y$ are not fibres of $V$.
In this case, we can extend the Seifert fibring
of $V$ to a Seifert fibring of $V \cup Y$.  Either $V \cup Y = S^3$, or
$V \cup Y$ has boundary.
\begin{itemize}
\item If $V \cup Y = S^3$ then we know by the classification of Seifert fibrings of $S^3$ that any fibring of $S^3$ has at most two singular fibres.  If $V$ is the complement of
a regular fibre of a Seifert fibring of $S^3$, then $V$ is a torus knot complement
$V \simeq M(0,1;\frac{\alpha_1}{\beta_1},\frac{\alpha_2}{\beta_2})$ with 
$\alpha_1\beta_2 - \alpha_2\beta_1 = \pm 1$. Otherwise, $V$ is the complement of a singular fibre, meaning that $V$ is a solid torus $M(0,1;\frac{\alpha}{\beta})$.
\item If $V \cup Y$ has boundary, we can repeat the above argument. 
Either $V \cup Y$ is a solid torus, in which case 
$V \simeq M(0,2;\frac{\alpha_1}{\beta_1})$
or a component of $S^3 \setminus \overline{V \cup Y}$ is a solid torus, so 
we obtain $V$ from the above manifolds by removing a Seifert fibre. By induction,
we obtain $V$ from either a Seifert fibring of a solid torus, or a Seifert 
fibring of $S^3$ by removing fibres. Thus either
$V \simeq M(0,n;\frac{\alpha_1}{\beta_1},\frac{\alpha_2}{\beta_2})$
for $n \geq 1$ with $\alpha_1\beta_2 -\alpha_2\beta_1 = \pm 1$,
$V \simeq M(0,n;\frac{\alpha_1}{\beta_1})$ for $n \geq 1$, 
or $V \simeq M(0,n;)$ for $n \geq 2$.  In order, these are the cases where
we remove only regular fibres from a fibring of $S^3$, regular fibres from
a fibring of a solid torus, and regular fibres plus the singular fibre
from a fibring of a solid torus.
\end{itemize}
\end{enumerate}
\end{proof}
\end{prop}

\begin{prop}\label{seifert-glt}
Each link in $S^3$ whose complement admits a Seifert-fibring 
is isotopic to some $S(p,q|X)$, excepting only the
key-chain links.

Given $p,q \in \Zed \setminus \{0\}$ let 
$p'=p/GCD(p,q)$, $q'=q/GCD(p,q)$
and let $m,l \in \Zed$ satisfy $p'm-lq'=1$, then
the complement of $S(p,q|X)$ is diffeomorphic to:
\begin{itemize}
\item $M(0,GCD(p,q);\frac{m}{q'},\frac{l}{p'})$ provided
$X=\emptyset$.
\item $M(0,1+GCD(p,q);\frac{m}{q'})$ provided $X=\{*_1\}$.
\item $M(0,1+GCD(p,q);\frac{l}{p'})$ provided $X=\{*_2\}$
\item $M(0,2+GCD(p,q);)$ provided $X=\{*_1,*_2\}$.
\end{itemize}
The the complement of the key-chain link $H^p$ is diffeomorphic
to $M(0,p+1;)$.

Let $A$ denote an index-set for the components of $S(p,q|X)$
which are neither $*_1$ nor $*_2$. Then the strong Brunnian
property of $S(p,q|X)$ is given by:
\begin{itemize}
\item $\{ \{*_1\} \}$ for the unknot $S(1,1|)$
\item $\{ \{*_1\}, \{*_2\} \}$ for the Hopf link $H^1 = S(2,2|)$
\item $\{ \{a\} : a \in A \} \cup X'$ if $p | q$ or $q | p$ but $p \neq 0$ and $q \neq 0$, where $X'$ is the collection of singleton subsets of $X$.
\item $X'$ if $p\nmid q$ and $q\nmid p$ 
where $X'$ is the collection of singleton subsets of $X$.
\end{itemize}
\begin{proof} 
Except for the Brunnian properties, this result follows immediately
from Proposition \ref{sf2}. To see the Brunnian properties, observe that
the linking number between regular fibres of $S(p,q|X)$ is $p'q'$,
the linking number between $*_1$ and $*_2$ is one, and the linking number
between a regular fibre and either $*_1$ or $*_2$ is $p'$ and $q'$ 
respectively.
\end{proof}
\end{prop}

Proposition \ref{seifert-glt} first appears in the literature in the paper of
Burde and Murasugi \cite{BM}, where they classify links in $S^3$ whose compliments
are Seifert-fibred. 

\begin{cor}\label{fibreslope}
Let $L$ be a link in $S^3$ such that $C_L$ admits a Seifert fibring. 
Provided $L$ is not the unknot nor the Hopf link, the fibring
is unique up to isotopy.
\begin{proof}
$S(p,q|X)$ is the unknot if and only if $X=\emptyset$, 
$GCD(p,q)=1$ and either $p=\pm 1$ or $q=\pm 1$.
$S(p,q)|X)$ is the Hopf link if and only if
either $X=\emptyset$ and $p=\pm q = \pm 2$ or
$|X|=1$ with $p=\pm q = \pm 1$. 
Thus, the complements of $S(p,q|X)$ which are
not unknot complements or Hopf link complements
all have the form
$M(0,b;\frac{\alpha_1}{\beta_1},\frac{\alpha_2}{\beta_2})$
where the sum of the number of boundary components plus the
number of singular fibres is at least $3$ with
$\alpha_1\beta_2-\alpha_2\beta_1=\pm 1$.
That these manifolds have unique Seifert-fibrings up to orientation-preserving
diffeomorphism follows from Theorem 2.3 in Hatcher's notes \cite{Hatcher3}.
Consider a horizontal essential annulus $S$ in one of these manifolds.
We have the basic relation among Euler characteristics 
$\chi(B) - \frac{\chi(S)}{n} = \sum_i(1-\frac{1}{\beta_i})$
where $B$ is the base space. Since $\chi(S)=0$, this equation has
no solution.  Thus by Proposition 1.11 of \cite{Hatcher3} all essential
annuli are vertical, therefore any diffeomorphism of these manifolds
is fibre-preserving.
\end{proof}
\end{cor}

\begin{prop}\label{Spquniq}
Let $\theta$ be the unique involution of the set $\{*_1,*_2\}$. 
Let $\mathcal{C} = \{S(p,q|X) : (p,q) \in \Zed^2, X \subset \{*_1,*_2\}\}$.
The equivalence relation $\sim$ of unoriented 
isotopy on $\mathcal{C}$ is generated by the relations: 
\begin{enumerate}
\item $S(p,q|X)\sim S(-p,-q|X) \ \ \forall (p,q)\in (\Zed \setminus \{0\})^2, X\subset \{*_1,*_2\}$.
\item $S(p,q|X) \sim S(q,p|\theta(X)) \ \ \forall (p,q)\in (\Zed \setminus \{0\})^2, X\subset \{*_1,*_2\}$.
\item $S(p,q|X \cup \{*_1\}) \sim S(p+\frac{p}{q},q+1|X\setminus \{*_1\})
\ \ \forall q | p$ with $q>0$, $X \subset \{*_1,*_2\}$.
\item $S(p,q|X) \sim S(-p,q|X)$ if either (a) 
$X=\emptyset$, $GCD(p,q)=1$ with $p = \pm 1$ or $q= \pm 1$
(b) $X=\emptyset$, $p, q = \pm 2$ or
(c) $|X|=1$, $\pm p = \pm q = 1$. This is the case where
$S(p,q|X)$ is an unknot or Hopf link.
\end{enumerate}
\begin{proof}
Restrict to the sub-class of $\mathcal{C}$ consisting of
$S(p,q|X)$ which are not unknot nor Hopf links. These
complements have a unique Seifert-fibring by Proposition 
\ref{fibreslope}. If $X=\emptyset$, $p\nmid q$ and $q\nmid p$, 
observe that items (1) and (2) generate $\sim$, this is 
because the classification of Seifert-fibred spaces 
up to fibre-equivalence  (Proposition 2.1 in \cite{Hatcher3}) 
tells us that $\sim$ is equivalent to the 
fibre-equivalence relation, thus we have proven 
more: up to isotopy there is only one orientation-preserving 
embedding of the complement of $S(p,q|\emptyset)$ in $S^3$.
If we broaden the class to include $X$ non-empty,
(1) and (2) still suffice to generate $\sim$ essentially
because $S(p,q|\emptyset)$ is contained as a sublink.
Consider a general $S(p,q|X)$. Let
$p'=p/GCD(p,q)$, and $q'=q/GCD(p,q)$. We know
$p', q'$ and $p'q'$ are the possible linking numbers of
the components of $S(p,q|X)$ thus we can determine whether
or not $p = \pm q$ via linking numbers.  For such a link,
the relative sign $p/q \in \{\pm 1\}$ can be computed
by coherently orienting three strands and computing the
linking number of two of them. Thus such links
are classified by the number of their components
together with the sign $p/q \in \{\pm 1\}$
which is equivalent to relations (1), (2) and (3). 
Consider the case $q|p$ but $p\nmid q$. 
In this
case $*_2 \in S(p,q|X)$ if and only some linking
number is $\pm 1$, and the linking number with $*_2$
would be $q'$. Thus each such is equivalent via
relations (1) through (3) to some
$S(q,p'q|X)$ where $X$ is either empty or contains $*_2$,
$q>0$ and $|p'|>1$. $q$ is the number of components of 
$S(q,p'q|X)$. Since we have assumed $S(q,p'q|X)$ is not the
unknot nor the Hopf link, $q+|X| \geq 3$ thus we can compute
$p'$ as a linking number of two coherently oriented strands
of $S(q,p'q|X)$ either in the complement of $*_2$ if $q=2$
or in the complement of another strand, thus relations
(1) through (3) suffice.

In the exceptional case of the unknot or the Hopf
link, $p/q \in \{\pm 1\}$ is not an invariant, thus
the exceptional relation (4).
\end{proof}
\end{prop}

\begin{defn}\label{fibreslopedef}
Given a Seifert-fibred $3$-manifold 
$M \subset S^3$ with $T \subset \partial M$
a boundary component, the fibre-slope of $M$ at $T$ is
$\frac{\alpha}{\beta}$ provided
$\pm(\alpha c_e + \beta c_i) \in H_1(T;\Zed)$
represents the homology class of a Seifert
fibre, where $c_e, c_i \subset T$ are external
and internal peripheral curves for $T$ such that
$lk(c_e',c_i)=+1$ where $c_e' \subset int(M)$
is parallel to $c_e$.

Given a link $L$ indexed by a set $A$ with
$C_L$ Seifert-fibred, for $a\in A$ define the
$a$-th fibre-slope of $L$ to be the fibre-slope
of $C_L$ at the boundary torus corresponding to $L_a$.
\end{defn}

The proof of the following is immediate.

\begin{prop}\label{FSEP}
If $M \subset S^3$ is Seifert-fibred with $T \subset \partial M$
a boundary torus and $f : M \to S^3$ the 
untwisted re-embedding of $M$, 
then the fibre-slope of $M$ at $T$ is the fibre-slope
of $f(M)$ at $f(T)$. 
Moreover, $M \cup N \subset S^3$ is Seifert-fibred
with $M \cap N = T$, then the fibre-slope of $M$ at $T$
is the inverse of the fibre-slope of $N$ at $T$.
\end{prop}

The importance of Proposition \ref{FSEP}
is that it gives us an obstruction to two Seifert-fibred manifolds
being adjacent in a link complement split along its
JSJ-tori.

\begin{prop}\label{FSCOMP}
The fibre-slopes of $S(p,q|X)$ are $LCM(p,q)/GCD(p,q)$ for
regular fibres, $p/q$ for $*_1$ and $q/p$ for $*_2$.
\end{prop}

The Thurston geometries \cite{ThuBook} on the Seifert-fibred
submanifolds of $S^3$ turn out to be non-unique. We give a sketch
of how to construct them. Milnor \cite{MilSP} has shown that 
any Seifert-fibred link complement in $S^3$ is the total space of a fibre bundle
over $S^1$.  
If $M$ is Seifert-fibred and $f : M \to S^1$ is a fibre-bundle,
let $F=f^{-1}(1)$ be the fibre.  The monodromy (attaching map),
as an element of the mapping-class group of $F$, is
always of finite order provided $F$ is not a torus.
This is because $F$ is essential in $M$ so we can
assume that either $F$ is a union of Seifert-fibres, or it
is transverse everywhere to the Seifert-fibres \cite{Hatcher3}.
If $F$ is transverse to the Seifert-fibres then
as an element of the mapping class group of $F$ ($\pi_0 \Diff(F)$) 
the monodromy has order equal 
to $|F \cap c|$ provided $c$ is a regular fibre.
If $F$ is a union of fibres, it is either $S^1 \times [0,1]$ or
a torus, and the mapping class group of $S^1 \times [0,1]$ is
finite. The interior of $\mathcal{S}_{g,b}$ is hyperbolic for
$2g+b \geq 3$, thus the monodromy is an isometry of the fibre for a suitable
hyperbolic metric on the fibre \cite{Ker}. Of
course, the only link complements that fibre over $S^1$ with
$F \simeq \mathcal{S}_{g,b}$ satisfying $2g+b<3$ are the
Hopf link and the unknot.
Theorem 4.7.10 of \cite{ThuBook} tells us that on top of having
the above finite-volume {$\bf H^2 \times E^1$}-structure, 
Seifert-fibred link complements also have an 
$\widetilde{PSL}(2,\Real)$-structure. 

Here is an example of how one can find the
{$\bf H^2 \times E^1$}-structures on the complement of
a torus knot. Let $C_{p,q} = S^3 \setminus T^{(p,q)}$. Think
of $T^{(p,q)}$ as the roots of the polynomial $f(z_1,z_2)=z_1^p-z_2^q$
where $(z_1,z_2)\in S^3 \subset \Complex^2$. The fibring
$g : C_{p,q} \to S^1$ is given by $g(z_1,z_2)=\frac{f(z_1,z_2)}{|f(z_1,z_2)|}$.
Define $C'_{p,q} = \Complex^2 \setminus f^{-1}(0)$.  There is
an action of the positive reals $\Real^+$ on $C'_{p,q}$ given by
$t.(z_1,z_2)=(t^{\frac{1}{p}}z_1,t^{\frac{1}{q}}z_2)$. 
As a function of $t$, $|t.(z_1,z_2)|$ is strictly increasing, thus 
$C'_{p,q} \simeq \Real^+ \times C_{p,q}$.
The function $f : C'_{p,q} \to \Complex \setminus \{0\}$ is a submersion
thus a locally-trivial fibre bundle.  If we let 
$C''_{p,q} = f^{-1}(S^1) \subset C'_{p,q}$ then 
$\Real^+ \times C''_{p,q} \simeq C'_{p,q}$ since $\Complex \setminus \{0\} \simeq \Real^+ \times S^1$.
So $C''_{p,q}$ and $C_{p,q}$ are homotopy-equivalent. 
Moreover, define $g' : C'_{p,q} \to S^1$
and $g'' : C''_{p,q} \to S^1$ by $g'(z_1,z_2)=\frac{f(z_1,z_2)}{|f(z_1,z_2)|}$
with $g''$ the restriction of $g'$.  This makes
the homotopy-equivalence $C''_{p,q} \to C_{p,q}$ a fibre 
homotopy-equivalence from $g''$ to $g$. Since both surfaces are $1$-ended,
they are diffeomorphic, moreover the homotopy-equivalence
preserves the peripheral structures of the fibres, so $g$ and $g''$
are smoothly-equivalent bundles. The fibre of $g''$, $g''^{-1}(1)$ is
$\{(z_1,z_2) \in \Complex^2 : z_1^p-z_2^q = 1 \}$ so the projection
map $\pi : \Complex^2 \to \Complex$ given by $\pi(z_1,z_2)=z_2$ 
restricts to a $p$-sheeted branched covering space with
$g''^{-1}(1)$ as the total space, $\Complex$ the base space,
having $q$ branch points in $\Complex$, all with ramification number $p$.
This makes $g''^{-1}(1)$ a surface of genus $(p-1)(q-1)/2$. 
The monodromy $g''^{-1}(1) \to g''^{-1}(1)$ can be written explicitly
as the map $(z_1,z_2) \longmapsto (e^{2\pi i/p}z_1,e^{2\pi i/q}z_2)$. 
As a covering transformation of $g''^{-1}(1)$ it 
is of order $pq$, with $p+q$ points of ramification, 
$p$ of which have ramification number $q$,
$q$ having ramification number $p$, branched over a disc with two
marked points. The idea for this computation comes from an 
example in an appendix of Paul Norbury's in the notes of 
Walter Neumann \cite{NeuNot}.

Other than Seifert-fibred manifolds, a primary source for
atoroidal manifolds is hyperbolic manifolds. Given an incompressible
torus $T$ in a complete hyperbolic $3$-manifold $M$, then
one would have an injection
$$\pi_1 T \to \pi_1 M \subset \Isom(\mathbf{H^3})$$
where $\Isom(\mathbf{H^3})$ is the group of isometries of
hyperbolic 3-space.
By the classification of hyperbolic isometries (see for example
Proposition 2.5.17 in Thurston's book \cite{ThuBook}) any
subgroup of $\Isom(\mathbf{H^3})$ isomorphic to $\Zed^2$
consists entirely of parabolic elements. The Margulis
Lemma, when applied to $\pi_1 T \to \Isom(\mathbf{H^3})$
(see for example \cite{ThuBook2} or \cite{Kap}) tells us that $M | T$
consists of two manifolds, one of which is diffeomorphic
to $T \times [0,\infty)$, thus $M$ is atoroidal.
Thurston went on to prove a rather sharp converse
\cite{Thurston}: the interior of a compact $3$-manifold $M$ with
non-empty boundary admits a complete hyperbolic metric if and
only if $M$ is prime, `homotopically atoroidal' and not the
orientable $I$-bundle over a Klein bottle, moreover the metric
is of finite volume if and only if $\partial M$ is a disjoint
union of tori.  `Homotopically atoroidal' means that subgroups
of $\pi_1 M$ isomorphic to $\Zed^2$ are peripheral. A convenient
refinement is that if $M$ is a compact, irreducible, atoroidal,
$3$-manifold with incompressible boundary
containing no essential annulus with $T \subset \partial M$ the
torus boundary components of $M$, then $M \setminus T$ admits a 
unique complete hyperbolic metric of finite volume with totally-geodesic
boundary \cite{Kap, Bon, ThuBook2}. If $\partial M = \emptyset$ then
$M$ is also known to be hyperbolic provided $M$ contains an 
incompressible surface.  Thus, hyperbolic geometry is a topological
invariant and it can be used to distinguish isotopy-classes of knots
and links. 
There is a corresponding theory of orbifolds that allows one to
go further and get sharp geometric invariants of 
smoothly-embedded graphs in $S^3$ \cite{Heard}.

\section{Companionship graphs for knots and links, splicing}\label{MAIN_RESULTS}

In this section we define two graphs associated to a link. The first,
called the `JSJ-graph,' (Definition \ref{def1}) describes the basic structure of the JSJ-decomposition of a link's complement.
We decorate the vertices of the JSJ-graph to get the `companionship graph' 
of a link (Definition \ref{COMP_GRAPH_DEF}), which, in Proposition \ref{IG_L_is_complete} we show is a complete isotopy invariant of links.
We turn our attention first to knots.  Companionship graphs for knots
have the particularly simple form of a rooted tree with vertices
labelled by `knot generating links' (Definition \ref{KGL_def}).
We develop a notion of splicing for knots
in Definition \ref{splicedef}. This allows us to construct knots 
with prescribed companionship graphs.
We describe the basic combinatorics of companionship graphs
for knots under splicing in Proposition \ref{KGL_EXCL}.
Theorem \ref{mainthm} gives a complete characterisation of companionship
graphs for knots, after which we give various examples.

We then turn our attention to companionship graphs of links.
We define the notion of a `splice diagram' in Definition \ref{splice_diagram}.
A splice diagram is not a concept of any real importance of its own, as it simply 
codifies some of the most elementary properties of a companionship graph.
The `Local Brunnian Property' (Proposition \ref{brun_excl}) is the first
fundamental property of companionship graphs. We show in Proposition
\ref{REAL_VALID} that splice diagrams satisfying the Local Brunnian Property
essentially encode for collections of disjoint embedded tori in link complements.
We proceed to call these diagrams `valid.' 
In Definition \ref{splice-relations} we give a revised notion of splicing
suitable for links.  The `fibre-slope exclusion property' (Lemma \ref{EXCEPTIONAL})
is the other main property satisfied by companionship graphs, as it 
encodes the minimality condition of the JSJ-decomposition.
Proposition \ref{graph_real} gives a complete characterisation of 
graphs that arrise as companionship graphs of links: they are the valid splice 
diagrams, labelled by Seifert-fibred and hyperbolic links satisfying the 
fibre-slope exclusion property. In Proposition \ref{graph_of_a_splice} we describe how companionship graphs behave under splicing, the most complicated case being
the case of splicing with the unknot.

\begin{defn}\label{def1} Given a topological space $X$, if
$\sim$ denotes the equivalence relation
$x \sim y \Longleftrightarrow$ there is a path from $x$ to $y$, 
define $[x]=\{ y \in X : y \sim x \}$ and $\pi_0 X = \{ [x] : x \in X\}$.

Given a non-split link $L$ in $S^3$ indexed by a set
$A$, let $T$ be the JSJ-decomposition of $C_L$,
indexed by a set $B$ disjoint from $A$.
The graph $G_L$ is defined to have vertex-set
$\pi_0 (C_L | T)$ and edge set $\pi_0 T$. 
We give $G_L$ the structure of a partially-directed graph, in that
some edges will have orientation. Given $b \in B$
let $M$ and $N$ be the two components of $C_L | T$
containing $T_b$. If $T_b$ bounds a solid torus $W$ 
in $S^3$ on only one side, and if $M \subset W$, then
we orient the edge $\pi_0 T_b$ so that its terminal
point is $\pi_0 M$.

For a split link $L$ with index set $A$, partition
$A$ as $A= A_1 \sqcup A_2 \sqcup \cdots \sqcup A_j$ so
that $C_L = C_{L_{A_1}} \# \cdots \# C_{L_{A_j}}$
is the prime decomposition of $C_L$. Define
$G_L$ for a split link $L$ to be $\sqcup_{i=1}^j G_{L_{A_i}}$.
$G_L$ is called the JSJ-graph of $L$.
\end{defn}

%%%%%%%%%%%%%%%%%%%% SIMPLE G_K example %%%%%%%%%%%%%%%%%%%%%%%%%

\begin{eg}\label{ex1} An example of a knot $K$ and its
JSJ-decomposition $T=\{T_1,T_2,T_3,T_4\}$.
$\K$ is a connected sum of a trefoil, 
a figure-8 and the Whitehead double of a figure-8 knot.  
{
\psfrag{k1}[tl][tl][0.8][0]{$T_1$}
\psfrag{k2}[tl][tl][0.8][0]{$T_2$}
\psfrag{k3}[tl][tl][0.8][0]{$T_3$}
\psfrag{k4}[tl][tl][0.8][0]{$T_4$}
\psfrag{K}[tl][tl][1][0]{$K$}
$$\includegraphics[width=15cm]{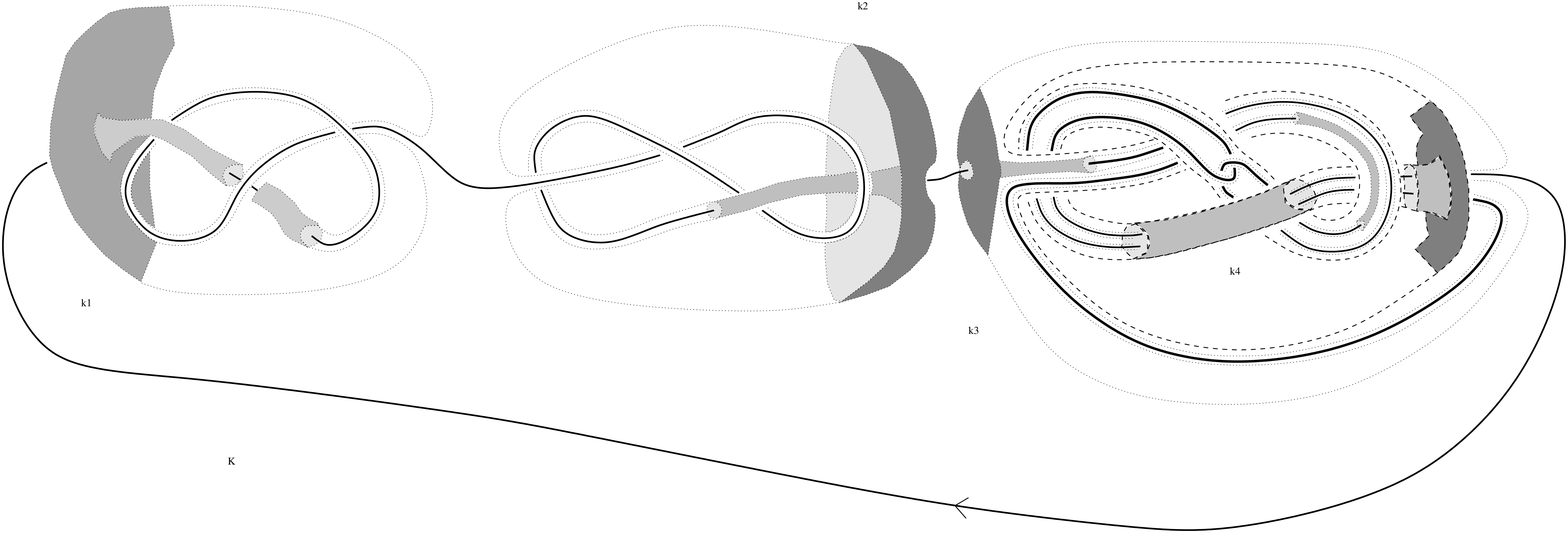}$$
}
Notice that the component of $C_K | T$ containing
$\partial C_K$ has as its untwisted re-embedding the 
complement of $H^3$. 
$T_1$ bounds a trefoil complement, 
$T_2$ a figure-8 complement, 
$T_3$ the Whitehead double of the figure-8 knot and 
$T_3 \cup T_4$ bounds a manifold, which, when re-embedded
is the complement of the Whitehead link. 
The Whitehead link complement and figure-8 complement are atoroidal since
they are finite-volume hyperbolic manifolds as mentioned at the
end of Section \ref{SEIFERT}.
{
\psfrag{t3}[tl][tl][0.7][0]{$\pi_0 T_1$}
\psfrag{t2}[tl][tl][0.7][0]{$\pi_0 T_2$}
\psfrag{t1}[tl][tl][0.7][0]{$\pi_0 T_3$}
\psfrag{tp}[tl][tl][0.7][0]{$\pi_0 T_4$}
\psfrag{jsj}[tl][tl][1][0]{$G_K$}
$$\includegraphics[width=7cm]{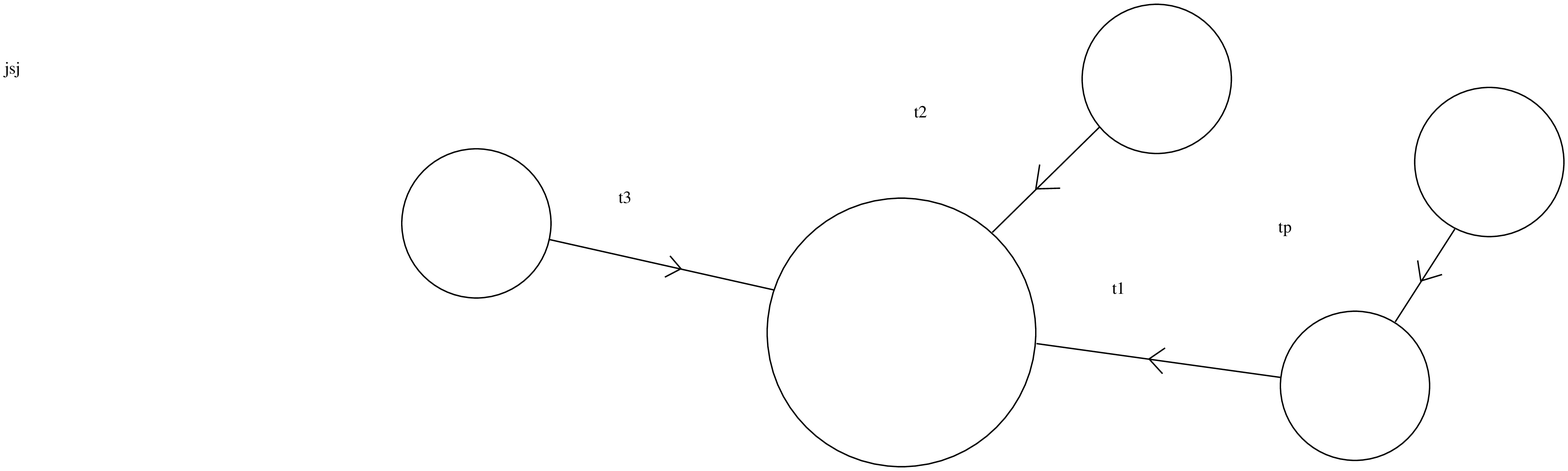} $$
}
\end{eg}

Graph constructions from the JSJ-decompositions of 3-manifolds were 
perhaps first made by Siebenmann \cite{Sieb}
in the context of homology spheres. Eisenbud and Neuman \cite{EN} 
made similar constructions for knots in homology spheres whose complements
are graph-manifolds.
%%%%%%%%%%%%%%%%%%% Definition IG_L

\begin{defn}\label{COMP_GRAPH_DEF}
Given a non-split link $L$ with index-set $A$, if $T$ is the
JSJ-decomposition of $C_L$, for each $v \in G_L$ let
$M(v) \in C_L | T$ be the component corresponding to $v$.
We define $\IG_L$ to be the partially-directed graph such that
each vertex $v$ is labelled by a link $\IG_L(v)$ satisfying:
\begin{enumerate}
\item If we forget the vertex labelling, $\IG_L$ is the JSJ-graph $G_L$.
\item The unoriented isotopy class of $\IG_L(v)$ is the companion
link to the component $M(v)$.
\item If $A(v)$ denotes the subset of $A$ corresponding to the components
of $M(v) \cap \partial C_L$, and $E(v)$ the subset of 
edges of $G_L$ incident to $v$, then $\IG_L(v)$ naturally indexed
by the set $A(v) \sqcup E(v)$. 
\item Given $v_1$ and $v_2$ adjacent vertices of $\IG_L$, let $\{e\}=E(v_1)\cap E(v_2)$,
then the orientation of $\IG_L(v_1)_e$ and $\IG_L(v_2)_e$ is chosen so that
if $f_1$ and $f_2$ are the untwisted re-embeddings 
$f_i : M(v_i) \to C_{\IG_L(v_i)}$, then
$lk(f_1^{-1}(l_1),\widetilde{f_2^{-1}(l_2)})=+1$ where 
$l_i \subset \partial C_{\IG_L(v_i)}$ is the standard longitude corresponding
to $\IG_L(v_i)_e$ respectively, and $\widetilde{f_2^{-1}(l_2)} \subset int(M(v_1))$
is a parallel translate of $f_2^{-1}(l_2)$.
\end{enumerate}
If $L$ is a split link, define
$\IG_L = \sqcup_{i=1}^k \IG_{L_{A_i}}$ where $C_L \simeq
C_{L_{A_1}} \# \cdots \# C_{L_{A_k}}$ is the prime decomposition
of $C_L$.
\end{defn}

Thus the union of the index sets for the links $\{\IG_L(v) : v \in \IG_L\}$ 
consist of $A$ together with the edges of $G_L$.  Counting with
multiplicity, exactly one component of the links $\{\IG_L(v) : v \in \IG_L\}$
is indexed by any element of $A$, and precisely two (corresponding to adjacent
vertices) are labelled by an edge of $G_L$.
Components decorated by elements of $A$ are called `externally-labelled.'
Components labelled by the edge-set of $G_L$ are `internally labelled.'
If $K$ is a knot, we consider it to be a link with index-set $\{*\}$.

Given an edge $e$ of $\IG_L$ with endpoints $v_1$ and $v_2$,
if we reverse both of the orientations of $\IG_L(v_1)_e$ and
$\IG_L(v_2)_e$, this would also satisfy the above definition, 
and $\IG_L$ is well defined modulo this choice.
%%%%%%%%%%%%%%%%%%%%% Souped-up of above example
{
\psfrag{1}[tl][tl][0.7][0]{$1$}
\psfrag{2}[tl][tl][0.7][0]{$2$}
\psfrag{3}[tl][tl][0.7][0]{$3$}
\psfrag{4}[tl][tl][0.7][0]{$4$}
\psfrag{S}[tl][tl][0.7][0]{$*$}
\psfrag{jsj}[tl][tl][1][0]{$\IG_K$}
$$\includegraphics[width=8cm]{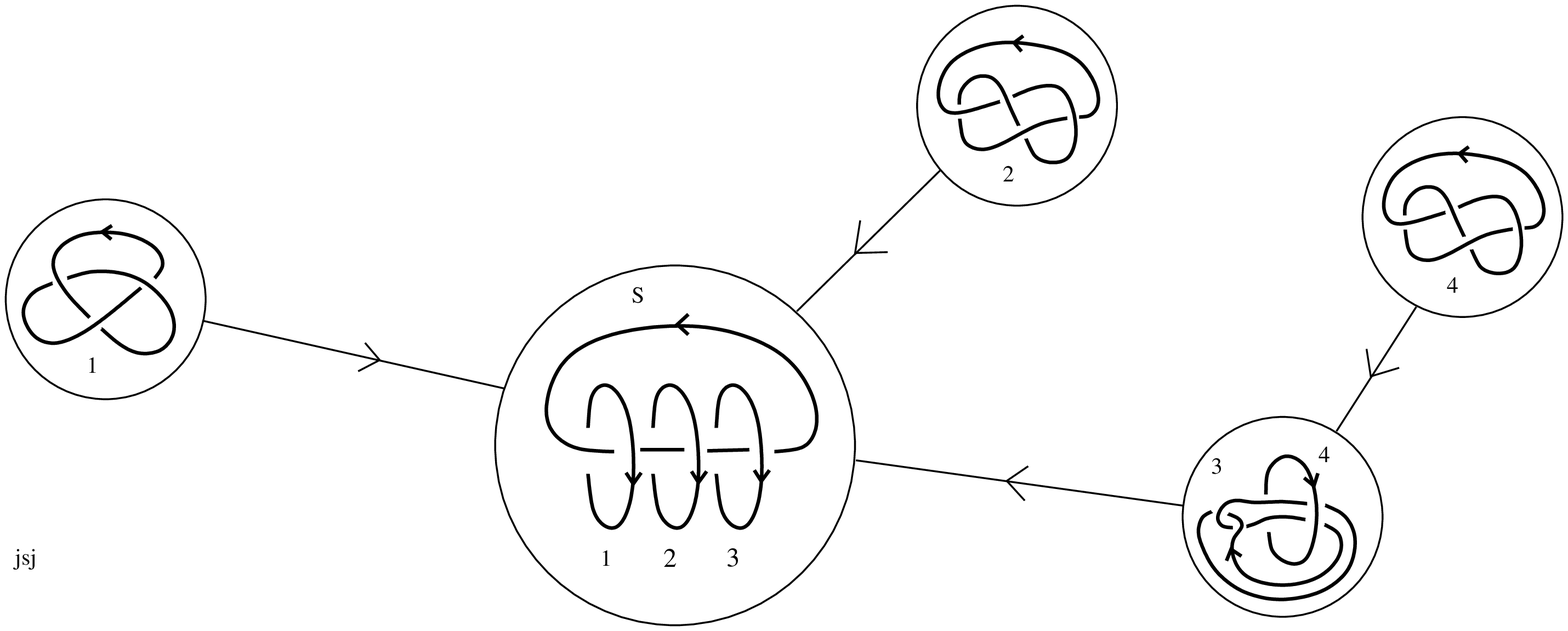} $$
}
%%%%%%%%%%%%%%%%%%% BASIC PROPS of IG_K, G_K, etc...
\begin{defn}\label{KGL_def}
A KGL indexed by $(A,b)$ is a link $L$ with index-set
$A \sqcup \{b\}$ such that $L_A$ is an unlink.
When it does not cause confusion, we will frequently
let $b=0$ and $A=\{1,2,\cdots,n\}$.
\end{defn}

Some elementary observations one can make about $\IG_L$:

\begin{itemize}
\item $\IG_L$ is connected if and only if $L$ is not split.
\item $\IG_L$ is acyclic -- ie: each component is a tree. 
This follows from the Generalised Jordan Curve Theorem \cite{Pollack},
since an embedded torus in $S^3$ separates.
\item If $K$ is a knot $\IG_K$ is a rooted tree, since
only one component of $C_K | T$ contains $\partial C_K$. 
By Proposition \ref{embtorbou}, for each vertex $v \in \IG_K$, 
$\IG_K(v)$ is a KGL. If we let $\IG_K(v)$ be indexed by $(A,b)$,
then the edges of $\IG_K$ corresponding to $A$ terminate at $v$. 
Provided $v$ is not the root of $\IG_K$, $b$ corresponds to an 
edge that starts at $v$ and terminates at its parent.
Thus, all the edges of $\IG_K$ are oriented and all sufficiently-long
directed paths in $\IG_K$ terminate at the root. 
\end{itemize}

\begin{defn}\label{splice_diagram}
Given a finite set $A$, a splicing diagram with external labels
$A$ is an acyclic, partially-directed graph $\IG$ such
that each vertex $v \in \IG$ is labelled by a link $\IG(v)$ whose
index-set is a subset of $\text{ (the edge-set of } \IG) \cup A$.
We demand that if $e \in \IG$ is an edge with $v_1,v_2 \in \IG$ its
endpoints, then one component of both
$\IG(v_1)$ and $\IG(v_2)$ is labelled by $e$, these are called
`internally-labelled'. We demand that for each $a \in A$ there exists a 
unique $v \in \IG$ such that $\IG(v)$ has a component indexed by $a$.
We denote this vertex by $v_a$, and we say $\IG(v_a)_a$ is
`externally-labelled'.

Given splicing diagrams $\IG$ and $\IG'$ with external index-sets
$A$ and $A'$ respectively, we say $\IG$ and $\IG'$ are equivalent
($\IG \sim \IG'$) if $A=A'$ and there exists
an isomorphism of partially-directed graphs $g : \IG \to \IG'$
together with unoriented isotopies $f(v)$ from 
$\IG(v)$ to $\IG'(g(v))$ for all $v \in \IG$ 
such that:
 \begin{itemize}
 \item If $e \in \IG$, $e' = g(e) \in \IG$, with 
$v_1,v_2 \in \IG$ the endpoints of $e$, and $v_i'=g(v_i)$ then
$f(v_i)(\IG(v_i)_e) = \epsilon \IG'(v_i')_{e'}$ for some
$\epsilon \in \{\pm\}$ valid for both $i \in \{1,2\}$. 
\item Given $a \in A$ and a component $\IG(v)_a$, then $f(v)(\IG(v)_a)=\IG'(g(v))_a$.
\end{itemize}
\end{defn}

See the conventions regarding unoriented isotopies in Section \ref{INTRODUCTION}
to make sense of the above definition of equivalence of splice diagrams.

\begin{prop}\label{IG_L_is_complete}
Two links $L$ and $Y$ are isotopic if and only if $\IG_L \sim \IG_Y$.
\begin{proof}
`$\Longrightarrow$' Is immediate since $\sim$ is an equivalence relation.

`$\Longleftarrow$' Let $h : S^3 \to S^3$ be
an isotopy from $L$ to $Y$. Thus $A=A'$ and $h(L_a)=Y_a$ for all $a \in A$,
moreover $h(C_L)=C_Y$.  Since the JSJ-decomposition is unique up to 
isotopy, then if $T$ is the JSJ-decomposition of $C_L$, we can assume
$h(T) \subset C_Y$ is the JSJ-decomposition of $C_Y$.  If we let
$M \in C_L | T$ then $h(M) \in C_Y | h(T)$ is isotopic to $M$, thus
by Proposition \ref{desplice_conv}, the companion link of 
$C_L$ corresponding to $M$ is
unoriented isotopic to the companion link of $C_Y | h(T)$, thus
we can let $f$ be this isotopy.
\end{proof}
\end{prop}

%%%%%%%%%%%%%%%%%%% SPLICING

The notion of `splicing' was first described by 
Siebenmann \cite{Sieb} in his work on JSJ-decompositions 
of homology spheres.  It was later adapted to the context of
links in homology spheres by Eisenbud and Neumann \cite{EN}.  
We give a further refinement of splicing, adapted specifically 
so that it constructs knots in $S^3$.

\begin{defn}\label{longdef}
A long knot is an embedding 
$f: \Real \times D^2 \to \Real \times D^2$ satisfying:
\begin{itemize}
\item $supp(f) \subset [-1,1] \times D^2$ i.e. $f$ is the identity outside of $[-1,1] \times D^2$.
\item The linking number of $f_{|\Real \times \{(0,0)\}}$ and $f_{|\Real \times \{(1,0)\}}$ is zero. 
\end{itemize}
\end{defn}

From a long knot $f$, one can construct a knot in $S^3$ in a canonical way. 
The image of $f_{|\Real \times \{(0,0)\}}$ is standard outside of 
$[-1,1]\times D^2$ so its one-point compactification is a knot in $S^3\equiv \dot{\Real^3}$.  This gives a bijective correspondence between isotopy
classes of long knots, and isotopy classes of knots.  The proof of this
appears in many places in the literature (\cite{BC, cubes} are recent examples)
and essentially amount to the observation that the unit tangent bundle
to $S^3$ is simply-connected. 

\begin{defn}\label{splicedef}
Let $\I = [-\infty,\infty]$, and let $L=(L_0,L_1,\cdots,L_n)$ an $(n+1)$-component KGL.  
Let $\tilde L$ be a closed tubular neighbourhood of $L$. 
Let $C_{L_i} = \overline{S^3 \setminus \tilde L_i}$ and
$C_u = \cap_{i=1}^n C_{L_i}$. Let $h=(h_1,h_2,\cdots,h_n)$ be a collection of disjoint 
orientation-preserving embeddings 
$h_i : \I \times D^2 \to C_{L_i}$ such that
$img(h_i) \cap \partial C_u = 
 img(h_{i|\I \times \partial D^2})$ with 
$h_i(\{0\}\times S^1)$ an oriented longitude for $L_i$.
We call $h$ a disc-system for $\tilde L$. 

Given $J=(J_1,J_2, \cdots, J_n)$ an $n$-tuple of non-trivial knots in $S^3$, 
 let $f=(f_1, f_2, \cdots, f_n)$ be their associated long knots. 
The re-embedding function associated to $L$, the disc-system $h$ and
knots $J$ is an embedding $R_h[L,J] : C_u \to S^3$ defined by:
$$R_h[L,J]= \left\{
\begin{array}{lr}
(h_i \circ f_i \circ h_i^{-1})(x) & \text{ if } x \in img(h_i) \\ 
x & \text{ if } x \in \overline{C_u \setminus \sqcup_{i=1}^n img(h_i)} \\
\end{array}
\right.$$
The splice of $J$ along $L$ is defined as
$$J \splice L = R_h[L,J](L_0).$$
\end{defn}

\begin{eg}\label{ex2}
{
\psfrag{l0}[tl][tl][1][0]{$L_0$}
\psfrag{l1}[tl][tl][1][0]{$L_1$}
\psfrag{l2}[tl][tl][1][0]{$L_2$}
\psfrag{h1}[tl][tl][1][0]{$h_1$}
\psfrag{h2}[tl][tl][1][0]{$h_2$}
\psfrag{l0ss}[tl][tl][1][0]{$L_0 \subset C_u$}
$$\includegraphics[width=7cm]{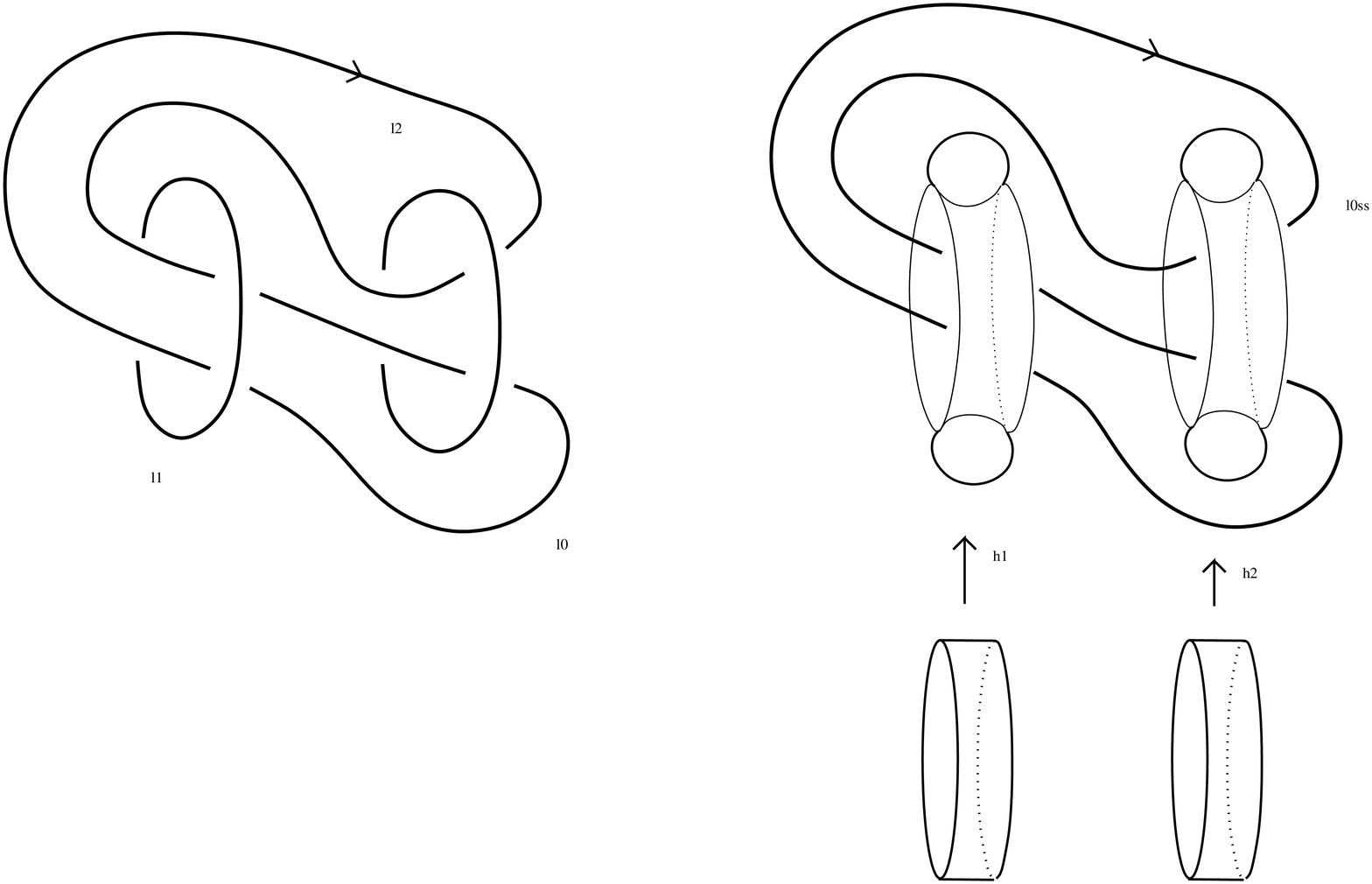}$$
}
If we use the figure-8 knot for 
$J_1$ and the knot $6_3$ from Rolfsen's knot 
table for $J_2$, then $K=J \splice L$ is the knot illustrated below.
{
\psfrag{K}[tl][tl][1][0]{$K = J \splice L$}
$$\includegraphics[width=4cm]{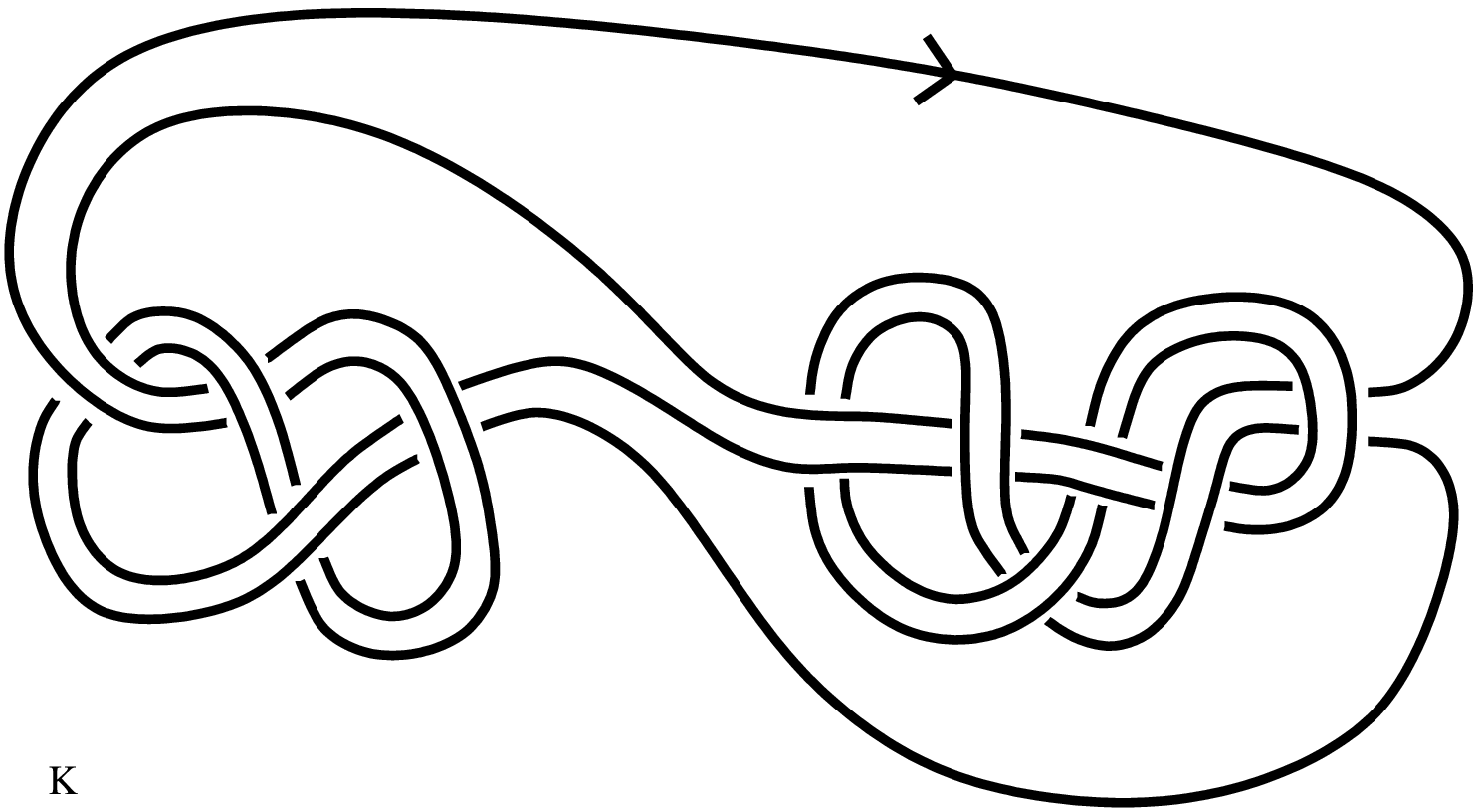}$$
}
\end{eg}

We will show that the isotopy-class of $R_h[L,J]$ does not depend
on $h$.  To do this, we show any two disc systems are related by a 
finite-sequence of  `elementary moves,' and that disc systems
related by a single elementary move give rise to isotopic
re-embedding functions.

\begin{defn}
Let $h$ and $h'$ be disc systems for $\tilde L$. 
An elementary move on disc $j$ from $h$ to $h'$ is a 
$1$-parameter family of embeddings
$H_i(t) : \I \times D^2 \to \overline{S^3 \setminus \tilde L_i}, 
\ i \in \{1,2,\cdots,n\}, \ t \in [0,1]$
such that:
\begin{itemize}
\item $H_i(0) = h_i$ and $H_i(1)=h'_i$ for all $i \in \{1,\cdots,n\}$.
\item $img(H_i(t)) \cap \partial \tilde L_i = 
img(H_i(t)_{|\I \times \partial D^2})$ for all 
$i \in \{1,\cdots,n\}, t \in [0,1]$.
\item $H_i(t) = H_i(0)$ for all $i \neq j$ and $t \in [0,1]$.
\end{itemize}
\end{defn}

\begin{prop}\label{re-embedding-lemma}
If two disc-systems $h$ and $h'$ are related by an elementary
move, then $R_h[L,J]$ is isotopic to $R_{h'}[L,J]$. 
\begin{proof}
Assume there is an elementary move on disc $j$ from $h$ to $h'$.
Extend $H_i(t) \circ f_i \circ H_i(t)^{-1}$ to the unique
embedding $\overline{S^3 \setminus \tilde L_i} \to S^3$ with support contained
in the image of $H_i(t)$. Let $\xi_i(t) : C_u \to S^3$ be its restriction.
Define $R_{H(t)}[L,J]$ by the formula:
$$R_{H(t)}[L,J] := 
\xi_{j} \circ 
\left( \xi_1(t) \circ \xi_2(t) \circ \cdots \circ \xi_{j-1}(t) \circ 
             \xi_{j+1}(t) \circ \cdots \circ \xi_{n-1}(t) \circ \xi_n(t)\right)
	       $$
This is well-defined and smooth since $\xi_i(t)(C_u)\subset C_u$ for all
$i \neq j$ and $t \in [0,1]$. By design $R_{H(0)}[L,J]=R_h[L,J]$ and
$R_{H(1)}[L,J]=R_{h'}[L,J]$.
\end{proof}
\end{prop}

The construction of the above isotopy is formally analogous
to the action of the operad of little $2$-cubes on the space of long 
knots \cite{cubes}.

Given an $(n+1)$-component KGL $L=(L_0,L_1, \cdots, L_n)$ we can 
represent a disc-system for $\tilde L$ by an embedding
$f : D \to S^3$ where $D=\sqcup_{i=1}^n D_i$ is a disjoint union
of $2$-discs and $\partial D_i = L_i$, as a closed regular 
neighbourhood of $D \cap C_u$ is such a disc-system.

\begin{prop}\label{re-embedding} 
Any two disc-systems $h$ and $h'$ for $\tilde L$ are related by
a sequence of elementary moves. 
\begin{proof}
Assume $f : D \to S^3$ is a disc system, intersecting $D$ 
transversely.  Consider the curves of intersection
$D \cap f(D) \subset D$.  If $D \cap f(D)=\partial D$, 
$f$ is isotopic to $D$ as a disc-system. 

So consider an innermost circle $C$ of 
$D \cap f(D)$ in $D$. 
$C$ is the boundary of two discs $C=\partial S_1$
and $C=\partial S_2$ where $S_1 \subset D$ and 
$S_2 \subset f(D)$.  
$S_1 \cup S_2$ is a sphere and so there is a unique 
$3$-ball $B \subset S^3$ with 
$\partial B = S_1 \cup S_2$ with 
$B \cap f(D) = S_2$. Let's assume $S_2 = B \cap f(D_i)$
Thus, there is an isotopy of
$f(D_i)$, supported in a neighbourhood of $B$ which
lowers the number of intersections of $f(D_i)$ with $D$.
This isotopy is a sequence of elementary moves on
$f$, one for every component of $img(f) \cap int(B)$.
By induction, we have the necessary collection of elementary
moves from $f$ to $D$.
\end{proof}
\end{prop}
\begin{defn}
Let $\mathbb G'$ be a sub-graph of $\IG_L$, thus it describes
some subset of $C_L | T$. 
If $M$ is the union of these submanifolds, and $f : M \to S^3$
the untwisted re-embedding, $f(M)$ is the complement of a tubular
neighbourhood some $1$-dimensional submanifold $X \subset S^3$.
Moreover, the boundary of $M$ are tori either of $C_L \cap M$, or they are 
edges $e \in \IG_L \setminus \mathbb G'$ incident to $\mathbb G'$. 
Define $\IG_L(\mathbb G')$ to be $X$, with components indexed by 
$A' \subset A$ corresponding to $C_L \cap M$, and the
edges of $\IG_L \setminus \mathbb G'$ incident to $\mathbb G'$.  Moreover, since
the external peripheral curves of these tori are naturally oriented
by the definition of $\IG_L$, this gives the components of 
$\IG_L(\mathbb G')$ a natural orientation.  We call $\IG_L(\mathbb G')$ the companion
link to $L$ for the sub-graph $\mathbb G' \subset \IG_L$.
\end{defn}

\begin{prop}
Given a knot $K$, let $L$ be the KGL decorating the root
of $\IG_K$, and let $\IG'_1, \cdots, \IG'_n$ be the sub-trees of
$\IG_K$ rooted at the children of $L$ in $\IG_K$, then
$$K \sim (\IG_K(\IG'_1),\cdots,\IG_K(\IG'_n))\splice L.$$
Moreover, given a vertex $v \in \IG_K$ with $\IG'$ the
maximal rooted sub-tree of $\IG_K$ rooted at $v$, then
$\IG_{\IG_K(\IG')} = \IG'$.
\end{prop}

Thus, the vertices of $\IG_K$ are an index-set for the
JSJ-companion knots to $K$, and any companion knot to
$K$ is a summand of a JSJ-companion knot.

We mention, without proof, a result on Alexander polynomials of
spliced knots.

\begin{thm} \cite{BZ} 
Given a knot $K$ let $\Delta_K(t)$ denote the Alexander polynomial of $K$.  
If $K= J \splice L$ where $L$ is a KGL, then 
$$\Delta_K(t) = \Delta_{L_0}(t) \cdot \prod_{i=1}^n
\Delta_{K_i}(t^{lk(L_0,L_i)}). $$
Here $lk(L_0,L_i)$ is the linking number between $L_0$ and $L_i$ where $L=(L_0,L_1,\cdots,L_n)$.
\end{thm}

The proof of the above theorem is an application 
of Proposition 8.23 of \cite{BZ}. 

A small observation on when splicing produces the unknot.

\begin{prop}\label{trivsplice} Let $L=(L_0,L_1,\cdots,L_n)$ be a KGL, and
$J=(J_1,\cdots,J_n)$ an $n$-tuple of non-trivial knots.
$J \splice L$ is the unknot if and only if $L$ is the unlink.
\begin{proof}
`$\Longrightarrow$' By design $\sqcup_{i=1}^n C_{J_i}$ naturally embeds in
$C_{J \splice L}$. As in the proof of Proposition \ref{embtorbou}, a
spanning-disc $D$ for $J \splice L$ can be isotoped off 
$\sqcup_{i=1}^n C_{J_i}$, giving a spanning disc $D'$ for 
$L_0$ disjoint from $\sqcup_{i=1}^n L_i$.
If $D'$ intersects the spanning discs for
$\sqcup_{i=1}^n L_i$ one can modify $D'$ through embedded surgeries
along the spanning discs of $\sqcup_{i=1}^n L_i$, resulting
in a spanning disc for $L_0$ disjoint from the spanning
discs for $\sqcup_{i=1}^n L_i$.

`$\Longleftarrow$' Let $D$ be a spanning disc for $L_0$ is disjoint
from $\sqcup_{i=1}^n L_i$. Then $R_h[L,J](D')$ is a spanning disc for
$J \splice L$.
\end{proof}
\end{prop}

%%%%%%%%%%%%%%%%%%%% GENERAL THEORY of IG_K

\begin{prop}\label{KGL_EXCL}Assume 
$L=(L_0,L_1,\cdots,L_n)$ is a non-compound KGL
and $J=(J_1,\cdots,J_n)$ is an $n$-tuple of non-trivial knots.
Provided both of the following statements are false:
\begin{itemize} 
\item $L$ is a Hopf link.
\item $L$ is a key-chain link and at least one of the
knots $J_1, \cdots, J_n$ is not prime.
\end{itemize}
then the root of $\IG_{J \splice L}$ is decorated by $L$, and the
maximal sub-trees of $\IG_{J \splice L}$ rooted at the children
of $L$ are $\IG_{J_1}, \cdots, \IG_{J_n}$ respectively.
\begin{proof}
First, consider the case that $L$ is not Seifert fibred. In this
case, the complement of $J \splice L$ is the union of
$R_h[J,L](C_L)$ and $C_{J_i}$ along the tori
$\partial C_{J_i}$ for $i\in \{1,2,\cdots,n\}$. $T_i$ is incompressible
in $C_{J_i}$ by the Loop Theorem. It is also incompressible in $R_h[J,L](C_L)$
since if it was not, an unknot would split off $L$ but we have assumed
$L$ is non-compound. Thus, the tori $T_i$ are incompressible in
$C_{J \splice L}$, moreover if we take the union of the collection
$\{T_1,\cdots,T_n\}$ together with the JSJ-decompositions of
$C_{J_i}$ for $i \in \{1,2,\cdots,n\}$ we get a collection of tori
$T$ such that $C_{J \splice L}|T$ consists of Seifert-fibred and 
atoroidal manifolds.  This is a minimal collection by assumption.

So we have reduced to the case $C_L$ Seifert-fibred. Provided
the roots of $\IG_{J_i}$ are all decorated by non-Seifert fibred
spaces, the above argument applies. So assume the root of
$\IG_{J_i}$ is also Seifert fibred. The fibre-slope of
$L_i$ is either $p/q$ if $L$ is $S^{(p,q)}$ or $0$ if
$L$ is a key-chain link. 
Consider the possible fibre-slopes of the relevant component
of the link decorating the root of $\IG_{J_i}$. It could either
be $\infty$ in the case of a key-chain link, or $LCM(r,s)/GCD(r,s)$
in the case of a torus knot or Seifert link. Thus, the only way the
Seifert-fibring could extend is the key-chain-key-chain case as in
all other cases, the fibre-slopes are not reciprocal (since $p\nmid q$
and $LCM(r,s)/GCD(r,s) \in \Natural$).
\end{proof}
\end{prop}

\begin{thm}\label{mainthm}
Given a knot $K$, the companionship tree is a connected splice 
diagram $\IG_K$ with external label $\{*\}$, such that every edge is
oriented, and every maximal directed path terminates at 
$v_*$, giving $\IG_K$ the structure of a rooted tree with root $v_*$.
\begin{enumerate}
\item Each vertex of $\IG_K$ is labelled by a link from the list:
 \begin{enumerate}
 \item Torus knots $T^{(p,q)}$ for $GCD(p,q)=1$, $p\nmid q$ and $q\nmid p$. 
 \item Seifert links $S^{(p,q)}$ for $GCD(p,q)=1$, $p\nmid q$. 
 \item Right-handed key-chain links $H^p$ for $p \geq 2$. 
 \item Hyperbolic KGLs.
 \item The unknot. 
 \end{enumerate}
\item Given a vertex $v \in \IG_K$, then $\IG_K(v)$ is some
KGL indexed by $(A,b)$ where the edges of $\IG_K$ corresponding
to $A$ are oriented towards $v$. If $v$ is not the root, then
$b$ corresponds to an edge oriented away from $v$.
\item If any vertex is decorated by a key-chain link $H^p$, 
none of its children are allowed to
be decorated by a key-chain link. 
\item A vertex of the tree $\IG_K$ is allowed to be decorated by the unknot 
if and only if the tree $\IG_K$ consists of only one vertex.
\end{enumerate}

The above properties are complete, in the sense that any graph 
satisfying the above properties is realisable as $\IG_K$ for some 
knot $K$. Moreover:

\begin{itemize}
\item Two knots $K$ and $K'$ are isotopic if and only if $\IG_K \sim \IG_{K'}$.
\item If $\IG_K$ consists of more than one vertex, then 
$K = \left(J_1,\cdots,J_n\right) \splice L$ where
the root of $\IG_K$ is labelled by $L$ and the knots $(J_1,\cdots,J_n)$
correspond to the maximal sub-trees rooted at the children of
the root of $\IG_K$.
\item There exists hyperbolic KGLs with arbitrarily many components. Thus one can 
realise any finite, rooted-tree as $G_K$ for some knot $K$ in $S^3$.  
\end{itemize}

\begin{proof} 
(1) Propositions \ref{seifert-glt} lists the Seifert-fibred links
and their Brunnian properties. The only Seifert-fibred KGL that 
we excluded from the list is the Hopf link. All remaining
non-compound KGLs are hyperbolic, by Thurston's 
Hyperbolisation Theorem \cite{Thurston}. 
(2) Follows from the definition.
(3) Is Proposition \ref{KGL_EXCL}.
(4) The unknot is the only knot whose complement does
not have an incompressible boundary.
Given a labelled rooted tree satisfying (1)--(4), 
one constructs the knot $K$ inductively on the height of
the tree, using Proposition \ref{KGL_EXCL} as the inductive step.

The first two bulleted ($\bullet$) points follow from
Proposition \ref{IG_L_is_complete}, Proposition \ref{embtorbou}
and Definition \ref{splicedef}.
The last bulleted point follows from a theorem of Kanenobu's \cite{Kan},
as will be explained. Consider a link $L$ indexed by a set $A$.  
The Brunnian property of $L$, $\Brun_L \subset 2^A$ satisfies 
the `Brunnian Condition':
\begin{itemize}
\item $S,T \in \Brun_L \text{ and } S \cap T \neq \emptyset \Longrightarrow S \cup T \in \Brun_L.$
\end{itemize}
For convienience, consider knots and the empty link to be split (this is a convention of Kanenobu):
\begin{itemize}
\item $\emptyset \notin \Brun_L$.
\item $\{a\} \notin \Brun_L \ \forall a \in A$. 
\end{itemize}
DeBrunner \cite{DeBrun} proved that if $\Brun \subset 2^A$ is any collection of subsets of a 
finite set $A$ satisfying the above three `Brunnian Conditions,' then there exists a link $L$ 
indexed by $A$ such that $\Brun_L = \Brun$.  Kanenobu went further \cite{Kan}: there exists
$L$ such that $\Brun_L = \Brun$ and for all $S \in \Brun$, 
$L_S$ is hyperbolic. Moreover, if $A \in \Brun$, then one can
assume all the components of $L$ are unknotted.

If we let $\Brun = \{ A \}$, then Kanenobu's
theorem gives us a hyperbolic KGL with $|A|$ components.
\end{proof}
\end{thm}

Theorem \ref{mainthm} allows one to consider the above class of rooted trees
as an index-set for the path-components of the space of embeddings
of $S^1$ in $S^3$, $\pi_0 \Emb(S^1,S^3)$. It is via this indexing that
the homotopy type of each component of $\Emb(S^1,S^3)$ is described in
the paper \cite{knotspace}.

Knots whose compliments have non-trivial JSJ-decompositions are 
quite common.  If one generates knots via random walks in $\Real^3$
where the direction vector is chosen by a Gaussian distribution,
then one typically gets a connect-sum \cite{Jun}.

The observation that if a knot in $S^3$ has a non-trivial JSJ-decomposition then it is
the splice of a link in $S^3$ with a knot in $S^3$
was first made in Proposition 2.1 of the Eisenbud-Neumann book \cite{EN}.  
Their point of view on the subject did not keep track of the Brunnian 
properties of the links involved, in that there is no analogue of 
Proposition \ref{unknottedcomp} in their work, as their work focuses
on links in homology spheres.
%%%%%%%%%%%%%%%%%%%%% EXAMPLES OF IG_K

\begin{cor}\label{basicsplices} Here are some elementary characterisations of 
some basic knot operations in terms of splicing.
\begin{itemize}
\item A knot $K$ is a connected-sum of $n$ non-trivial knots for $n \geq 2$
if and only if $K=J\splice L$ where $L$ is an $(n+1)$-component key-chain link, 
and the knots $J=(J_1,J_2,\cdots,J_n)$ are all non-trivial.
\item A knot $K$ is a cable knot if and only if $K$ is a splice knot 
$K=J\splice L$ such that $L$ is a $(p,q)$-Seifert link. We consider a torus
knot to be a cable of the unknot. Thus $K$ is a proper cable knot (a cable of a
nontrivial knot) if and only if $K=J \splice L$ where $L$ is a non-trivial knot.  
\item A knot $K$ is a (untwisted) Whitehead double if 
$K=J \splice L$ where $L$ is the Whitehead link.
\end{itemize}
\end{cor}
Corollary \ref{basicsplices} also appears in Schubert's work \cite{Sch2}. 

We give various examples of spliced knots and companionship trees.
Let $F_8$ denote the figure-8 knot.

{
\psfrag{cableknot}[tl][tl][0.9][0]{ $ K = T^{(-3,2)} \splice S^{(2,17)}$}
\psfrag{gk}[tl][tl][1][0]{$\IG_K$}
\psfrag{cs}[tl][tl][0.9][0]{$ K=((T^{(2,3)}, F_8) \splice H^2) 
\splice S^{(2,-17)} $}
\hskip 15mm
$\includegraphics[width=5cm]{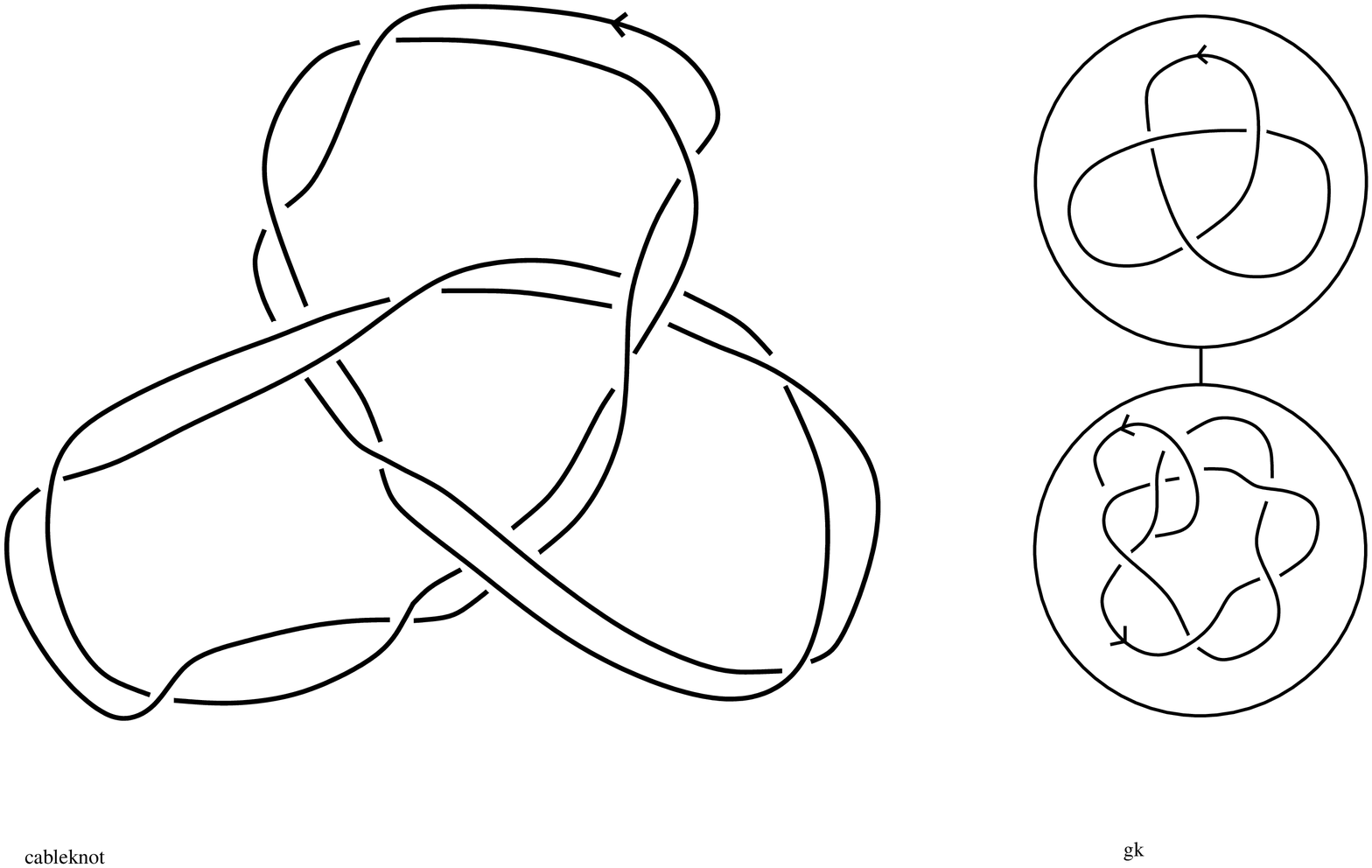}$
\hfill
$\includegraphics[width=6cm]{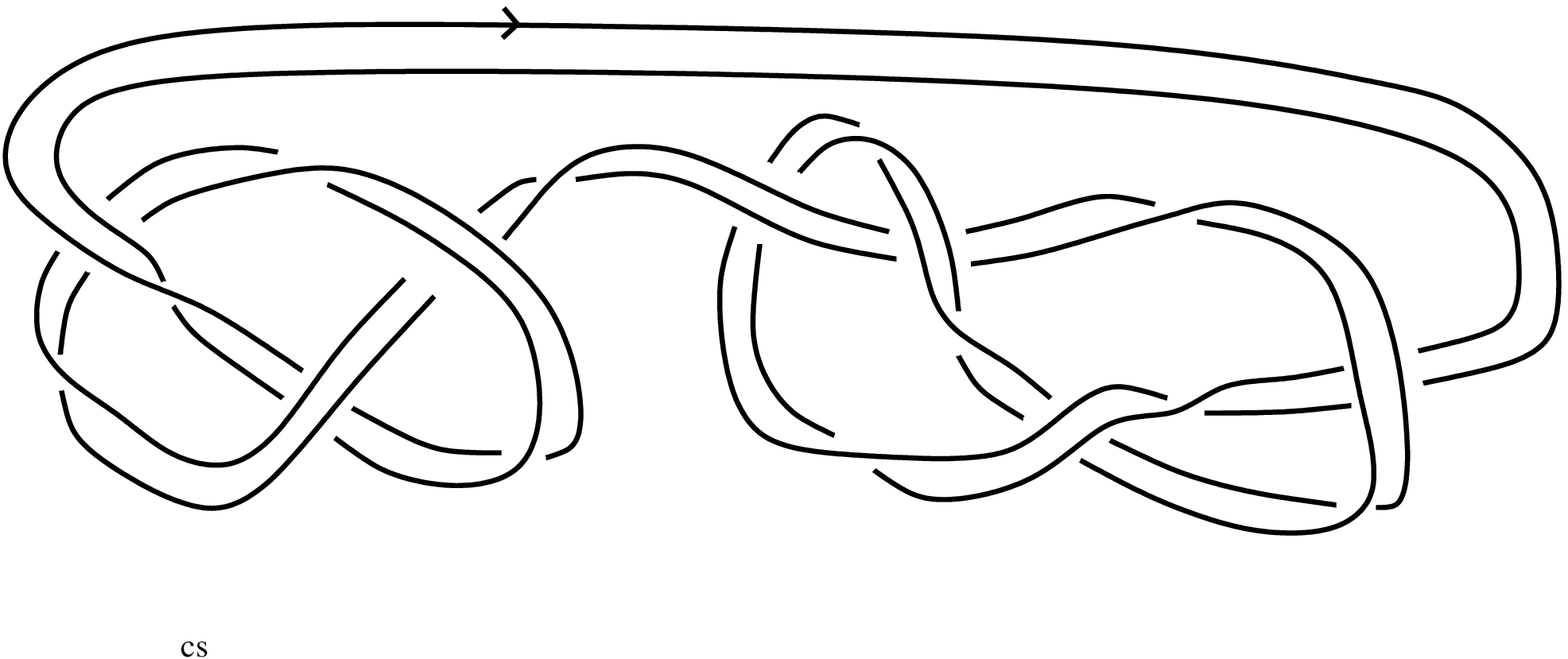}$
\hskip 15mm
}
\vskip 5mm
{

}
Let $W$ denote the Whitehead link.
Let $B=(B_0,B_1,B_2)$ denote the Borromean rings.
Let $B(i,j)$ be the $3$-component 
link in $S^3$ obtained from $B$ by
doing $i$ Dehn twists about the spanning disc of $B_1$ and $j$ Dehn twists
about the spanning disc for $B_2$. 

{
\psfrag{GHL}[tl][tl][1][0]{$\IG_K $}
\psfrag{n}[tl][tl][1][0]{$K=F_8 \splice W$}
\hskip 7mm
$\includegraphics[width=6cm]{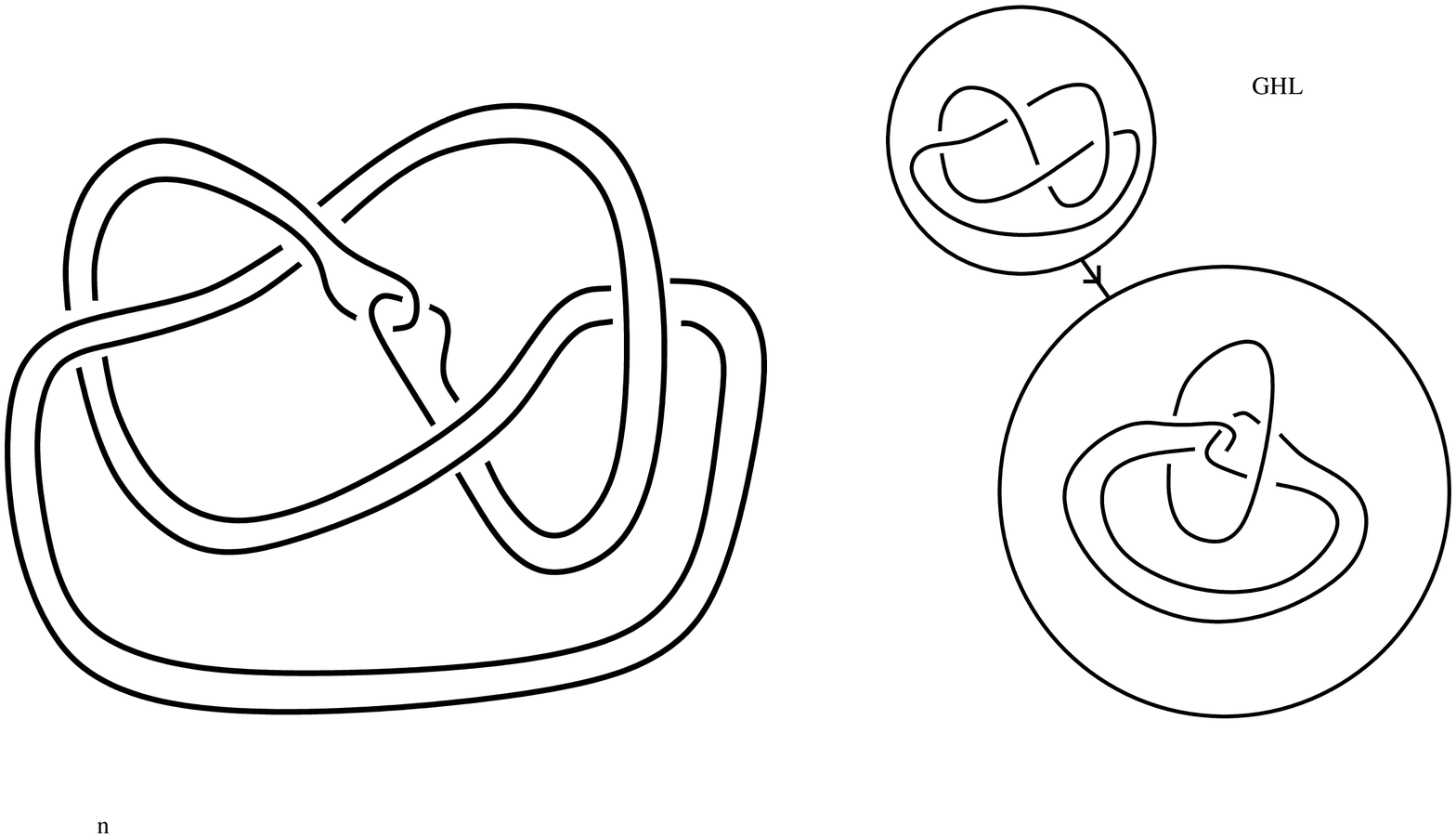}$
\hfill
\psfrag{n}[tl][tl][0.9][0]{$K=(F_8,T^{(3,2)}) \splice B_{0,6}$}
\psfrag{T}[tl][tl][1][0]{$\IG_K$}
$\includegraphics[width=7cm]{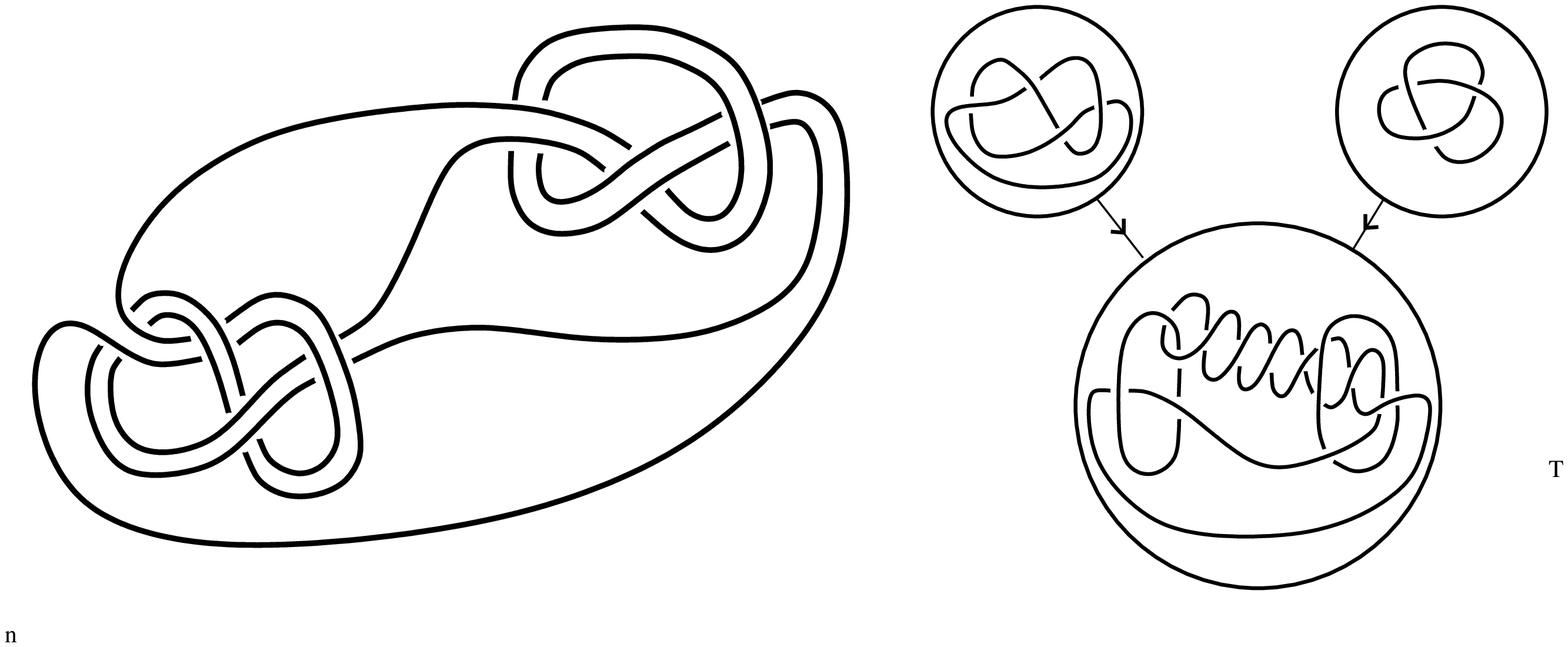}$
\hskip 7mm
}
\vskip 5mm

%%%%%%%%%%%%%%%%%%%%% BASIC PROPS OF IG_L, G_L

The graphs we associate to links in $S^3$
have more complicated combinatorics for three reasons:

\begin{itemize}
\item Link complements are not prime provided the link is split. 
This will result in our graphs being a union of disjoint trees.
\item There are link complements with incompressible tori
that separate components of the link, thus the associated graphs
are not always rooted.
\item The tori in the JSJ-decomposition of a link complement are
not always knotted. 
\end{itemize}

An example of a link with an unknotted torus in its JSJ-decomposition.
{ 
\psfrag{1}[tl][tl][0.8][0]{$1$}
\psfrag{2}[tl][tl][0.8][0]{$2$}
\psfrag{a}[tl][tl][0.8][0]{$a$}
\psfrag{3}[tl][tl][0.8][0]{$3$}
\psfrag{4}[tl][tl][0.8][0]{$4$}
\psfrag{L}[tl][tl][1][0]{$L$}
\psfrag{IGL}[tl][tl][1][0]{$\IG_L$}
$$\includegraphics[width=12cm]{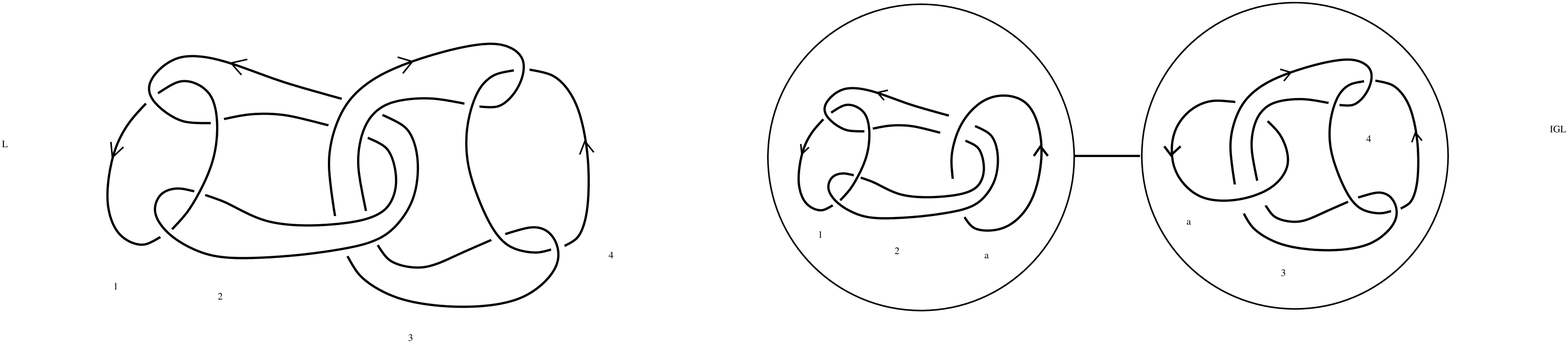}$$
}
%%%%%%%%%%%%%%%%%%%%% SPLICING IG_L's
\noindent The rest of this section will be devoted to describing the class of 
labelled graphs that can be realised as $\IG_L$ for a link $L$ in $S^3$,
and how the graphs behave under the corresponding notion of
splicing.
To do this, we identify the local rules that allow us to determine
if a link $Y$ can decorate a vertex in a graph $\IG_L$ for some $L$.

%%%%%%%%%%%%%%%%%%%%%% GENERAL THEORY IG_L, G_L

Given a vertex $v \in \IG_L$ in a companionship graph $\IG_L$
of a link $L$ with index set $A$, we partition the components of 
$\IG_L(v)$ into four classes:
\begin{enumerate}
\item Those which correspond to an
oriented edge of $\IG_L$ whose terminal point is $v$.
\item Those which correspond to an
oriented edge of $\IG_L$ whose initial point is $v$.
\item Those that correspond to an unoriented edge of 
$\IG_L$ incident to $v$.
\item The components $\IG_L(v)_a$ of $\IG_L(v)$ which 
have labels in the set $A$.
\end{enumerate}

Components of type (1) through (3) are indexed by a subset of $B$,
the index-set for the tori in the JSJ-decomposition of $C_L$.
Components of $\IG_L(v)$ of type (4) are indexed by a subset of $A$.

\begin{prop}{\it (Local Brunnian Exclusion Property)}\label{brun_excl}
Let $L$ be a link in $S^3$. Fix a vertex $v$ of $\IG_L$
and let $A(v)$ be the index-set for $\IG_L(v)$. 
Let $A_1, A_2, A_3 \subset A(v)$ be the indices corresponding
to items (1), (2) and (3) above. Then the following
statements hold:
\begin{enumerate}
\item $A_1 \in \BrunS_{\IG_L(v)}$.
\item $\forall a \in A_2, A_1 \cup \{a\} \notin \BrunS_{\IG_L(v)}$
\item $\forall a \in A_3, A_1 \cup \{a\} \in \BrunS_{\IG_L(v)}$
\end{enumerate}
\begin{proof}
Let $V$ be the submanifold of $C_L | T$ corresponding to
$v$. We index the tori of $\partial V$ by the set $A'$.
For each $a \in A_1$, $\partial_a V$ is a torus which
bounds a knotted solid torus $J_a$ in $S^3$ containing $V$.
The collection $\{J_a : a \in A_1\}$ have disjoint complements,
so by Proposition \ref{embtorbou}, $A_1 \in \BrunS_{\IG_L(v)}$.

Fix $a' \in A_2 \cup A_3$. $\partial_{a'} V$ bounds a solid torus
in $S^3$ disjoint from $V$. This solid torus is a neighbourhood
of some knot $K_{a'}$ in $S^3$. 
Provided $A_1 = \{a_1,a_2,\cdots,a_n\}$, define 
$L' = (\IG_L(v)_{a'}, \IG_L(v)_{a_1}, \cdots, \IG_L(v)_{a_n})$.
Then by the definition of our identification $V \simeq C_{\IG_L(v)}$,
$K_{a'} = J \splice L'$ where $J=(J_{a_1},\cdots,J_{a_n})$. 
By Proposition \ref{trivsplice} 
$K_{a'}$ is unknotted if and only if
$L'$ is an unlink, if and only if
$A_1 \cup \{a'\} \in \BrunS_{\IG_L(v)}$. This proves points (2) and (3).
\end{proof}
\end{prop}

A splice diagram will be called valid if it satisfies the 
Local Brunnian Property.  Valid splice diagrams essentially keep
track of embedded tori in link complements, as we will show
in the next proposition.

\begin{defn}
Given a link $L$ with index-set $A$ and a family of disjoint
embedded tori $T \subset C_L$ indexed by a set $B$, we define
the splice diagram associated to the pair $(L,T)$ to be the
graph whose underlying vertex-set is $\pi_0 (C_L | T)$, 
edge-set is $\pi_0 T \equiv B$, such that each vertex $v$ is decorated by
a link following Definition \ref{COMP_GRAPH_DEF} points (3) and (4).
We orient the edges of this graph following Definition \ref{def1}.
The notation we use for this graph is $\IG_{(L,T)}$.
\end{defn}

Notice that if $L$ is not a split link, then $\IG_L = \IG_{(L,T)}$
provided $T \subset C_L$ is the JSJ-decomposition.

\begin{defn}\label{splice-relations}
Given links $L$ and $L'$ with index sets 
$A$ and $A'$ such that 
$\{a\} \in \BrunS_L$ and $a' \in A'$. 
Let $f : S^3 \to S^3$ be
an orientation-preserving diffeomorphism 
such that $f(R[L_a,L'_a](C_{L_a}))$ is a closed
tubular neighbourhood of $L'_{a'}$ such that
an oriented longitude of $L'_{a'}$ corresponds
to an oriented meridian of $C_{L_a}$.
We define the splice of $L'$ and $L$ to be:
$$L' \underset{a' \ a}{\splice} L = 
L'_{A' \setminus \{a'\}} \cup f(L_{A \setminus \{a\}})$$
and we index this link by the set
$(A \setminus \{a\})\sqcup (A'\setminus \{a'\})$.
Provided $A \cap A' = \{a\}$ we simply denote the splice
by $L' \splice L$
\end{defn}

Note that if both $\{a\} \in \BrunS_L$ and $\{a'\} \in \BrunS_{L'}$,
then $L' \underset{a' \ a}{\splice} L$ and 
 $L \underset{a \ a'}{\splice} L'$ are isotopic, so without
any harm we can consider the splicing notation for links 
to be symmetric 
$L \underset{a \ a'}{\splice} L' = L' \underset{a' \ a}{\splice} L$. 

{ 
\psfrag{sp}[tl][tl][0.8][0]{$\splice$}
\psfrag{eq}[tl][tl][0.8][0]{$=$}
\psfrag{1}[tl][tl][0.8][0]{$1$}
\psfrag{2}[tl][tl][0.8][0]{$2$}
\psfrag{3}[tl][tl][0.8][0]{$3$}
\psfrag{4}[tl][tl][0.8][0]{$4$}
\psfrag{5}[tl][tl][0.8][0]{$5$}
$$\includegraphics[width=12cm]{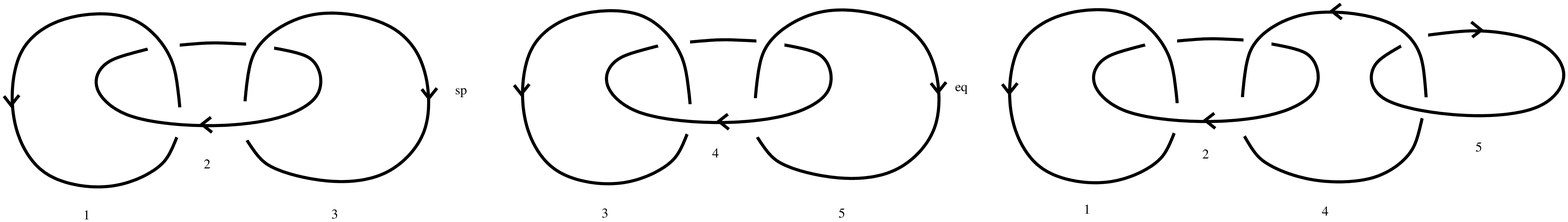}$$
}
\begin{prop}\label{REAL_VALID}Given a valid splice diagram $\IG$, then 
there exists link $L$ and tori $T \subset C_L$ such that
$\IG_{(L,T)} \sim \IG$, moreover, the pair $(L,T)$ is unique
up to isotopy in the sense that if $\IG_{(L',T')} \sim \IG$ then
there exists an isotopy $f$ from $L$ to $L'$ such that
$f(T)=T'$.
\begin{proof}
The proof of uniqueness is essentially the same as that of Proposition \ref{IG_L_is_complete}.
We prove existence by induction. For this we can assume $\IG$ is
connected, and since the initial step is true by design,
we proceed to the inductive step. 
Let $e$ be an edge of $\IG$ with
endpoints $v_1$ and $v_2$. $e$ partitions $\IG$ into two sub-graphs
$\IG'$ and $\IG''$.  

Provided $e$ is unoriented, then both $\IG'$
and $\IG''$ are valid splice diagrams therefore can be realised
by pairs $(L',T')$ and $(L'',T'')$ respectively.  

If $e$ is oriented, assume its terminal point is $v'$ and initial
point is $v''$. Assume $v' \in \IG'$ and $v'' \in \IG''$. Then
$\IG''$ is a valid splice diagram which can be realised
by a pair $(L'',T'')$.  $\IG'$ may not be a valid splice
diagram. Consider the sub-graph $Y$ of $\IG'$ defined recursively 
to be the sub-graph of $\IG'$ containing $v'$ such that if $w \in Y$
and if $f$ is an oriented edge of $\IG'$ whose initial point is in $Y$,
with $A_1 \cup \{f\} \in \BrunS_{\IG'(w)}$ then the endpoint of $f$
is also in $Y$.  If we unorient all the edges of $Y$ in $\IG'$, we
obtain a valid splice diagram which can be realised by some
pair $(L',T')$.  

By the definition of $L' \underset{e \ e}{\splice} L''$, 
$T'$ and $T''$ naturally correspond to a collection of
disjoint tori $T$ in the complement of $L' \underset{e \ e}{\splice} L''$.
Moreover, $\IG = \IG_{(L' \underset{e \ e}{\splice} L'',T)}$
\end{proof}
\end{prop}

\begin{defn}Given an oriented edge $e \in \IG$ in a valid splice diagram,
the sub-graph $Y$ of $\IG$ rooted at the terminal-point of $e$,
constructed in the proof of Proposition \ref{REAL_VALID} will be 
called the downward consequences of $e$.

Given two valid splice diagrams $\IG'$ and $\IG''$ with external
index-sets $A', A''$ realisable by
$(L',T')$ and $(L'',T'')$, provided $a' \in A'$ and $a'' \in A''$,
and either $\{a'\} \in \BrunS_{L'}$ or $\{a''\} \in \BrunS_{L''}$,
the splice diagram $\IG_{(L' \underset{a' \ a''}{\splice} L'',T)}$
constructed in Proposition \ref{REAL_VALID} will be denoted
$\IG' \BAR{a'}{a''} \IG''$.
\end{defn}

We note a convenient global property of companionship graphs,
allowing one to determine the Brunnian properties of a link via
that of its companions. It also shows how one can
determine the edge orientations of $\IG_L$ from the vertex-labels.

\begin{prop}\label{brun_of_a_splice}({\it Global Brunnian Property})
Given $L$ indexed by $A$ any edge $e \in \IG_L$ separates
$\IG_L$ into two sub-graphs $\IG'$ and $\IG''$. Let
$L' = \IG_L(\IG')$ and $L'' = \IG_L(\IG'')$ be the companions
indexed by the sets $A'$ and $A''$ respectively with $A'\cap A'' = \{e\}$
and $(A' \cup A'') \setminus \{e\} = A$.
Then $\BrunS_L$ is the set:
$$\{ B \subset A : 
(B \cap A' \cup \{e\} \in \BrunS_{L'} \text{ and } B \cap A'' \in \BrunS_{L''}) \text{ or }
(B \cap A'' \cup \{e\} \in \BrunS_{L''} \text{ and } B \cap A' \in \BrunS_{L'}\}$$
If $v' \in \IG'$ and $v'' \in \IG''$ are the endpoints of $e$, then
$e$ is oriented from $v'$ to $v''$ (resp. $v''$ to $v'$) if and only
if $\{e\} \notin \BrunS_{L'}$ and $\{e\} \in \BrunS_{L''}$ (resp.
$\{e\} \notin \BrunS_{L''}$ and $\{e\} \in \BrunS_{L'}$). Moreover,
$e$ is unoriented if and only if $\{e\} \in \BrunS_{L'} \cap \BrunS_{L''}$.
If $e$ is either unoriented or oriented from $v'$ to $v''$ (resp. $v''$ to $v'$)
then $\IG_{L'} \sim \IG'$ (resp. $\IG_{L''} \sim \IG''$). If $e$ is oriented
from $v'$ to $v''$ (resp. $v''$ to $v'$) then $\IG_{L''}$ (resp. $\IG_{L'}$)
is equivalent to $\IG''$ (resp. $\IG'$) after unorienting the edges of the
downward consequences of $e$.
\begin{proof}
Let $T \subset C_L$, $T' \subset C_{L'}$ and
$T'' \subset C_{L''}$ be the tori corresponding to $e$.
We prove the statement about $\BrunS_L$, the remaining statements
are corollaries of the proof of Proposition \ref{REAL_VALID}.

`$\subset$': Given $B \in \BrunS_L$, $T$ is compressible in $C_{L_B}$. 
By the Loop Theorem, either 
$T'$ is compressible in $C_{L'_{B \cap A' \cup \{e\}}}$ or 
$T''$ is compressible in $C_{L''_{B \cap A'' \cup \{e\}}}$.

If $T'$ is compressible in $C_{L'_{B \cap A' \cup \{e\}}}$
then since $L_B$ is an unlink one can choose the spanning discs
to be disjoint from $T$, thus
$L'_{B \cap A'}$ and $L''_{B \cap A''}$ are unlinks. 
$L'_a$ bounds a disc disjoint from $L'_{B \cap A'}$, therefore 
$L'_{B \cap A' \cup \{e\}}$ is also an unlink.
The argument for $T''$ compressible in $C_{L''_{B \cap A'' \cup \{e\}}}$
is formally identical as it is a symmetric argument.

`$\supset$': Let $\{M',M''\}=C_L | T$ with $f' : M' \to C_{L'}$ and
$f'' : M'' \to C_{L''}$ the untwisted re-embeddings.  
If $B \cap A' \cup \{e\} \in \BrunS_{L'}$ and $B \cap A'' \in \BrunS_{L''}$,
then if $D'$ and $D''$ are the spanning discs for $L'_{B \cap A'}$ and
$L''_{B \cap A''}$, $f'^{-1}(D') \cup f''^{-1}(D'')$ are spanning
discs for $L_B$. This is also a symmetric argument.
\end{proof}
\end{prop}

Notice that in our proof of the Global Brunnian Property, we did not
use the incompressibility of the tori $T$, thus it applies equally
well to valid splice diagrams.

\begin{eg}\label{ST1}
Consider the link below $L=(L_1,L_2)$.
{ 
\psfrag{L}[tl][tl][1][0]{$L$}
\psfrag{1}[tl][tl][0.8][0]{$1$}
\psfrag{2}[tl][tl][0.8][0]{$2$}
\psfrag{a}[tl][tl][0.8][0]{$a$}
\psfrag{b}[tl][tl][0.8][0]{$b$}
\psfrag{c}[tl][tl][0.8][0]{$c$}
\psfrag{d}[tl][tl][0.8][0]{$d$}
\psfrag{IGL}[tl][tl][1][0]{$\IG_L$}
\psfrag{e}[tl][tl][0.8][0]{$e$}
$$\includegraphics[width=14cm]{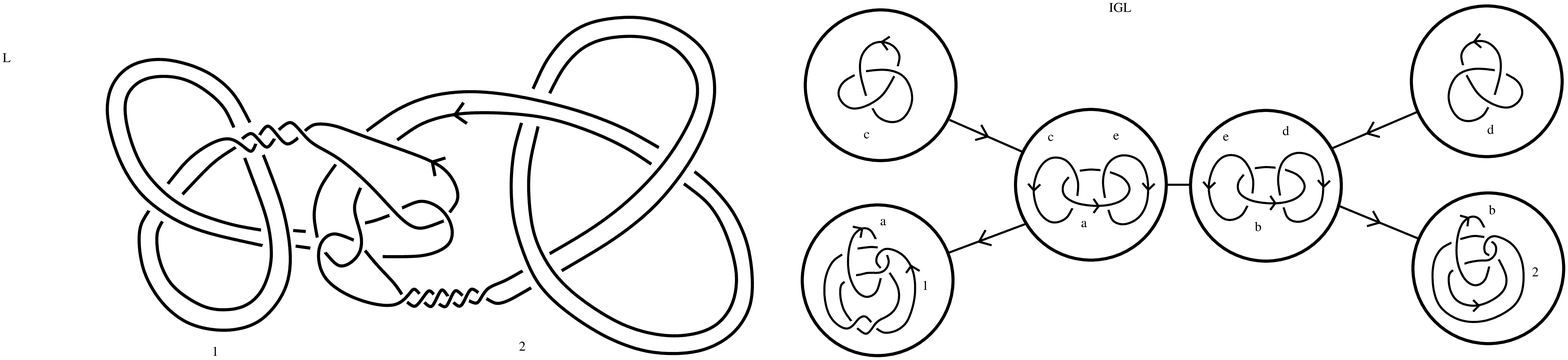}$$
}
Let $\IG'$ be the sub-graph of $\IG_L$ obtained by deleting the vertex
labelled by the (left-handed) trefoil indexed by $c$.
{ 
\psfrag{L}[tl][tl][1][0]{$L'=\IG_L(\IG')$}
\psfrag{1}[tl][tl][0.8][0]{$1$}
\psfrag{2}[tl][tl][0.8][0]{$2$}
\psfrag{a}[tl][tl][0.8][0]{$a$}
\psfrag{b}[tl][tl][0.8][0]{$b$}
\psfrag{c}[tl][tl][0.8][0]{$c$}
\psfrag{C}[tl][tl][1][0]{$c$}
\psfrag{IGL}[tl][tl][1][0]{$\IG_{L'}$}
\psfrag{d}[tl][tl][0.8][0]{$d$}
\psfrag{e}[tl][tl][0.8][0]{$e$}
$$\includegraphics[width=14cm]{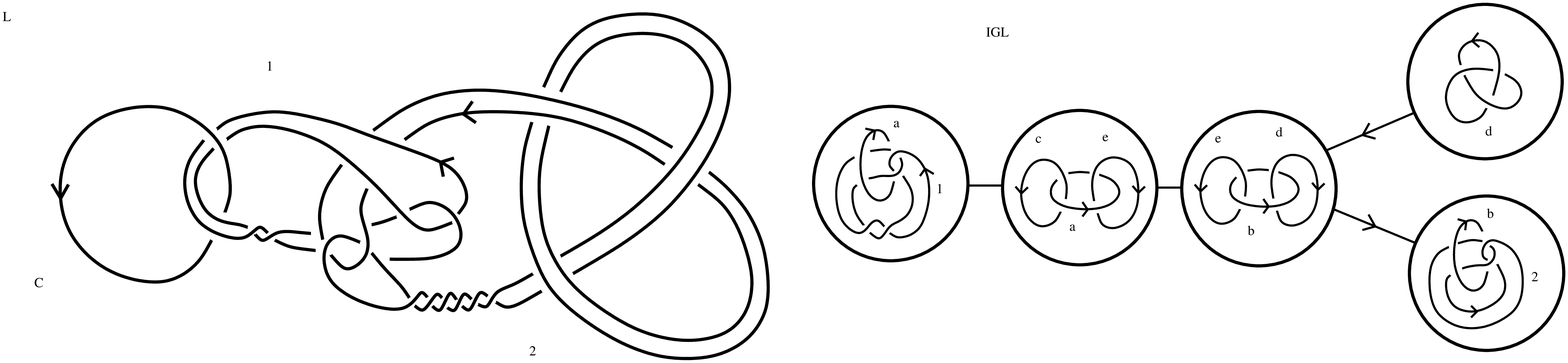}$$
}
\end{eg}

There is another local property satisfied by companionship
graphs. Only certain Seifert-fibred links may be adjacent
in $\IG_L$. As we have seen, given an edge $e \in \IG_L$ the
fibre-slopes of the components of the two adjacent links corresponding
to $e$ are never multiplicative inverses of each other (see 
Propositions \ref{FSEP}, \ref{KGL_EXCL}). 

\begin{lem}{\it (Fibre-slope Exclusion Property)}\label{FIB_EXCL}
A vertex decorated by an unknot or a Hopf link can not be
adjacent to any other vertex in a companionship graph $\IG_L$
for any link $L$.

Given two adjacent vertices in a companionship graph $\IG_L$ 
decorated by Seifert-fibred links whose unoriented isotopy
class representatives are $S(p,q|X)$ and $S(a,b|Z)$
respectively, there are several combinations that can not 
occur for any link $L$, which we list.

\begin{enumerate}
\item $S(p,q|\{*_1\})$ and $S(a,b|X)$ with edge corresponding
to $*_1 \in S(p,q|\{*_1\})$ and a regular fibre of $S(a,b|X)$,
provided $\frac{q}{p} = \frac{LCM(a,b)}{GCD(a,b)}$.
\item $S(p,q|X \cup \{*_1\})$ and $S(a,b|Z \cup \{*_2\})$
with edge corresponding to $*_1$ and $*_2$ respectively,
such that $\frac{p}{q}=\frac{a}{b}$.
\item $H^p$ and $H^q$ for $p,q \geq 2$ with edge corresponding 
to components of $H^p$ and $H^q$ having different fibre-slopes
ie: a `key' and a `keyring' component respectively.
\end{enumerate}
\end{lem}

The Hopf links and unknots are isolated since one can obtain
any rational number (or $\infty$) as the fibre-slope
of some Seifert fibring of the complement.  All other complements
have unique Seifert-fibrings (Proposition \ref{fibreslope}).

Lemma \ref{FIB_EXCL} can be recast into a splicing rule.

\begin{defn}\label{EXCEPTIONAL} Given links of the form $S(p,q|X)$ and $S(a,b|Z)$, 
let $r$ denote any regular fibre of $S(p,q|X)$ or $S(a,b|Z)$.
Let $L$ be a link with index-set $A$, and let $O$ be the
unknot with index set $\{*\}$. The following is the list of 
`exceptional splices':
\begin{enumerate}
\item  $S(a,b|X) \underset{r \ *_1}{\splice} S(p,q|\{*_1\}) = 
S\left(\frac{GCD(a,b)+GCD(p,q)-1}{GCD(a,b)}(a,b)|X\right)$ provided 
$\frac{q}{p} = \frac{LCM(a,b)}{GCD(a,b)}$.
\item $S(p,q|X \cup \{*_1\}) \underset{*_1 \ *_2}{\splice} S(a,b|Z \cup \{*_2\}) =
S\left(\frac{GCD(a,b)+GCD(p,q)}{GCD(a,b)}(a,b)|(X \setminus \{*_1\}) \cup (Z \setminus \{*_2\})\right)$ provided
$\frac{p}{q}=\frac{a}{b}$.
\item $H^p \underset{k \ g}{\splice} H^q = H^{p+q-1}$ provided
$k$ and $g$ represent a `key' and a `keyring' of $H^p$ and $H^q$
respectively.
\item $H^1 \splice L = L$ for any choice of components, provided
$H^1$ is the right-handed Hopf link.  If $H^1$ is left-handed, then
$H^1 \splice L$ consists is $L$ with one component's orientation
reversed.
\item $L \underset{a \ *}{\splice} O = L_{A \setminus \{a\}}$ for all
$a \in A$.
\end{enumerate}
\end{defn}

\begin{prop}\label{graph_real}Given a valid splice 
diagram $\IG$ satisfying the Fibre-slope exclusion
property such that each vertex $v \in \IG$ is labelled
by either a Seifert-fibred or hyperbolic link $\IG(v)$,
there exists a unique
link $L$ (up to isotopy) such that $\IG_L \sim \IG$.
Moreover, given any two links $L'$ and $L''$ with $a' \in A'$
and $a'' \in A''$ satisfying either $\{a'\} \in \BrunS_{L'}$
or $\{a''\} \in \BrunS_{L''}$, 
$\IG_{L' \underset{a' \ a''}{\splice} L''} = \IG_{L'} \BAR{a'}{a''} \IG_{L''}$ provided
$\IG_{L'} \BAR{a'}{a''} \IG_{L''}$ satisfies the
Fibre-Slope Exclusion Property.
\begin{proof}
This follows immediately from Proposition \ref{REAL_VALID}
and the proof of Proposition \ref{KGL_EXCL}.
\end{proof}
\end{prop}

We proceed to investigate the companionship tree
of $L' \underset{a' \ a''}{\splice} L''$ for arbitrary
links $L'$ and $L''$. To describe the case where one of
$L'$ or $L''$ is the unknot, we will need the operations
of `splitting' and `deletion' for valid splice diagrams,
defined below.

\begin{defn}
Consider a valid splice diagram $\IG$ with external
index-set $A$, vertex-set $V$ and edge-set $E$.

If $v \in V$ is such that $\IG(v)$ is split, we will define
a valid splice-diagram called the splitting of 
$\IG$ at $v$, denoted $\IG|v$.  
Let $X$ be the index-set for $\IG(v)$, and let
$X=\sqcup_{i=1}^k X_i$ be the partition of $X$ so that $C_{\IG(v)} \simeq 
C_{\IG(v)_{X_1}} \# \cdots \# C_{\IG(v)_{X_k}}$ is the prime decomposition
of $C_{\IG(v)}$.  
Define $\IG|v$ to be the splice diagram whose vertex-set is 
$(V \setminus \{v\}) \sqcup \left(\sqcup_{i=1}^k \pi_0 C_{\IG(v)_{X_i}}\right)$,
and whose edge set is $E$. We label the vertices 
$w \in V \setminus \{v\}$ by
$\IG(w)$, and the vertices $\pi_0 C_{\IG(v)_{X_i}}$ by
$\IG(v)_{X_i}$.

If $a \in A$, we define a valid splice-diagram
called the deletion of $a$ in $\IG$, denoted $\IG.a$.
Let $\IG'$ be the maximal sub-graph of $\IG$ with vertex-set $V \setminus v_a$.
If $A'$ is the index-set of $\IG(v_a)$, let $\IG''= \IG_{L_{A'\setminus \{a\}}}$.
Let $E'$ denote the edges of $\IG$ that are not edges of
$\IG'$, then $\IG' \sqcup \IG''$ is naturally a valid splice-diagram
once we append the edges $E'$.
\end{defn}

If one thinks in terms of the pair $(L,T)$ such that
$\IG_{(L,T)} = \IG$, splitting corresponds to finding a $2$-sphere
in the complement of $L \cup T$ that seperates components, while
$\IG.a$ corresponds to deletion of the component $L_a$ and
adding appending the tori from the JSJ-decomposition of $L_{A' \setminus \{a\}}$.

\begin{prop}\label{graph_of_a_splice}
Given two links $L'$ and $L''$ with index sets $A'$ and $A''$ respectively
such that $L' \underset{a' \ a''}{\splice} L''$ is defined, and
provided 
$\IG_{L'}(v_{a'}) \underset{a' \ a''}{\splice} \IG_{L''}(v_{a''})$ is
not an exceptional splice, 
$\IG_{L' \underset{a' \ a''}{\splice} L''}=\IG_{L'} \BAR{a'}{a''} \IG_{L''}$.

In the case of an exceptional splices (1) through (4), one obtains
$\IG_{L' \underset{a' \ a''}{\splice} L''}$ from
$\IG_{L'} \BAR{a'}{a''} \IG_{L''}$ by replacing the sub-graphs 
(see Definition \ref{EXCEPTIONAL}):
\begin{enumerate}
\item  $S(a,b|X) \text{ --- } S(p,q|\{*_1\})$ by
$S(\frac{GCD(a,b)+GCD(p,q)-1}{GCD(a,b)}(a,b)|X)$.
\item $S(p,q|X \cup \{*_1\}) \text{ --- } S(a,b|Z \cup \{*_2\})$
by $S(\frac{GCD(a,b)+GCD(p,q)}{GCD(a,b)}(a,b)|(X \setminus \{*_1\}) \cup (Z \setminus \{*_2\}))$.
\item $H^p \text{ --- } H^q$ by $H^{p+q-1}$.
\item $H^1 \text{ --- } L$ by $L$, or $L$ with the corresponding orientation reversed
if $H^1$ is the negatively oriented Hopf link.
\end{enumerate}
We call the operations (1) through (4), together with (5)
below `elementary reductions'.

In the case (5) of an exceptional splice with an unknot
$L' \text{ --- } \  O$, then $\IG_{L' \text{ --- }  O}$
is obtained from $\IG_{L'}.a$ by performing (recursively) 
all possible elementary reductions and splittings until the
there are no more available to do.
\begin{proof}
All cases except case (5) are direct consequences of Lemma \ref{FIB_EXCL}.
In case (5), $\IG_{L'}.a'$ describes a collection of tori $T$ dividing
the complement of $L'_{A' \setminus \{a'\}}$ into atoroidal and
Seifert-fibred spaces.  $C_{L'_{A' \setminus \{a'\}}}$ may not
be prime, and the tori may by compressible.  
If $S$ be a sphere
in $C_{L'_{A' \setminus \{a'\}}}$ intersecting $T$ transversely
separating components of $L'_{A' \setminus \{a'\}}$,
and $D$ is an innermost disc, then this disc is a spanning disc for
a component of one of the links decorating $\IG_{L'}.a'$, thus we
can perform a reduction of type (5). If there is no innermost disc,
we can perform a splitting. So after performing enough type (5) moves
and splittings, we are reduced to a disjoint union of diagrams that
describe incompressible tori in irreducible link complements that 
split the link complement into Seifert-fibred and atoroidal manifolds.
Any minimal such collection is the JSJ-decomposition, and 
moves (1) through (4) are by design what it takes to get to
that minimal collection.
\end{proof}
\end{prop}

Provided one never needs to use step (5), Proposition 
\ref{graph_of_a_splice} gives a rather
simple description of the companionship graph of a splice.
A good example of Proposition \ref{graph_of_a_splice} is
obtained by deleting a component of the link in
$L$ in Example \ref{ST1}.

In the case of a splice with an unknot $L' \text{ --- }  O$
where $L'$ is hyperbolic, Proposition \ref{graph_of_a_splice} 
says nothing. We mention a result 
of Thurston's that gives some insight into this situation.  
First, an example.

\begin{eg}\label{voleg} The link $L'$ contains the link $L$ from Example \ref{ST1}
as a sublink. 
{ 
\psfrag{L}[tl][tl][1][0]{$L'$}
\psfrag{2}[tl][tl][1][0]{$2$}
\psfrag{1}[tl][tl][1][0]{$1$}
\psfrag{3}[tl][tl][1][0]{$3$}
$$\includegraphics[width=6cm]{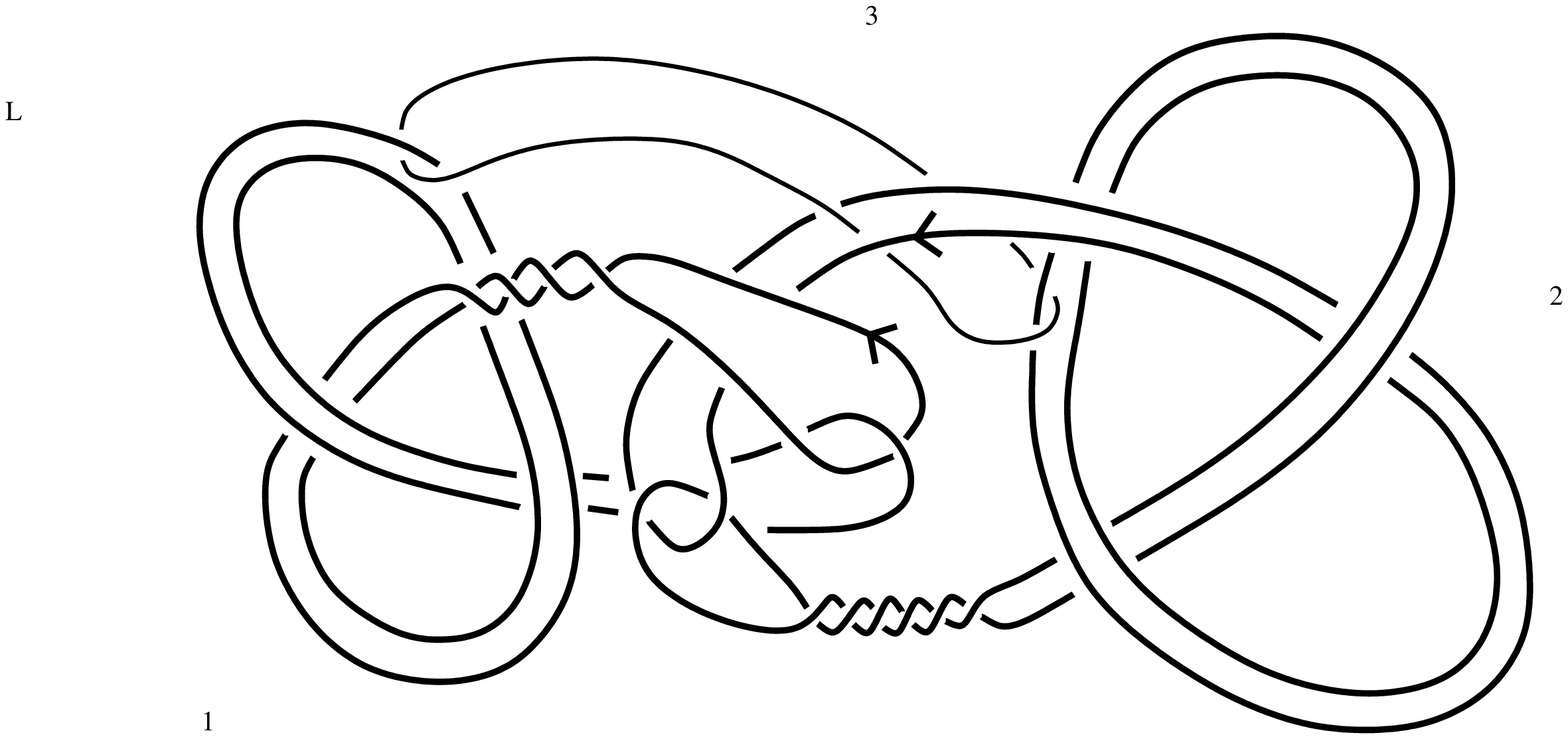}$$
}
$L'$ is hyperbolic with hyperbolic volume
$42.7594$. The link $L$ from Example \ref{ST1} has 
two hyperbolic companion links, both with hyperbolic
volume approximately $3.663$. 
\end{eg}

This is not a coincidence, as
we will explain.

\begin{defn}
Given a link $L$, the Gromov Norm of $L$ is defined to be the
sum 
$$|L| = \sum_{v \in \IG_L} 
\left\{
\begin{array}{lr}
Vol(S^3 \setminus \IG_L(v)) & \text{ if } S^3 \setminus \IG_L(v) \text{ is hyperbolic} \\
0 & \text{ if } S^3 \setminus \IG_L(v) \text{ is Seifert-fibred} \\
\end{array}
\right.$$
where $Vol(\cdot)$ is the hyperbolic volume.
\end{defn}

See \cite{Gromov} for a proper definition of Gromov Norm. That the 
proper definition reduces to the above in the case of link complements
is a theorem of Thurston's \cite{ThuBook2}. Thurston
also proved that if a link $L'$ is a proper sub-link of a
hyperbolic link $L$ then $|L'|<|L|$. 
Thus, if $L'$ is obtained
from an arbitrary link by deleting a component, $|L'|\leq |L|$
and one has equality if and only if when constructing $\IG_{L'}$
using the procedure in Proposition \ref{graph_of_a_splice},
one must never delete a component of a hyperbolic companion to $L$.

We end with one note on the computability of the companionship graph $\IG_L$. 
Typically one starts with a diagram for $L$ and then constructs a 
(topological) ideal triangulation of its complement, with an algorithm 
such as the one implemented in SnapPea \cite{Snap}.  The algorithm of Jaco, 
Letscher, Rubinstein \cite{JLR} computes the JSJ-graph $G_L$ from the
ideal triangulation. 
Link complements in $S^3$ are known to have solvable word
problems and these methods extend to an isotopy classification
of knots that seems (at present) to be not practical \cite{Wald}.
Commonly-used computer programs called Snappea \cite{Snap} and 
Orb \cite{Heard} 
frequently find the hyperbolic structures on hyperbolisable
$3$-manifolds (and orbifolds).
Once SnapPea has the hyperbolic metric, it then finds the canonical 
polyhedral decomposition \cite{Epstein} from which it can determine if
two hyperbolic links are isotopic by a simple combinatorial 
check. So in practise one can 
frequently use SnapPea to determine whether any two links are
isotopic. 
It's unfortunate that there is as of yet no formal justification for 
the effectiveness of these programs.
The Manning Algorithm \cite{Manning}  
takes a hyperbolisable 3-manifold with a solution to its word problem
and produces a hyperbolic metric, but this requires the usage of
solutions to word problems, which typically have very long run-times.

\providecommand{\bysame}{\leavevmode\hbox to3em{\hrulefill}\thinspace}

\Addresses

\end{document}